\documentclass{amsart}

\usepackage[utf8]{inputenc}

\usepackage{amssymb}
\usepackage{amsmath}
\usepackage{amsthm}
\usepackage{enumitem}
\usepackage{pgf,tikz}
\usetikzlibrary{arrows}

\numberwithin{equation}{section}

\theoremstyle{plain}
\newtheorem{thm}{Theorem}[section]
\newtheorem{lemma}[thm]{Lemma}
\newtheorem{prop}[thm]{Proposition}
\newtheorem{coroll}[thm]{Corollary}

\theoremstyle{definition}
\newtheorem{defn}[thm]{Definition}

\newtheorem{conj}[thm]{Conjecture}
\newtheorem{remark}[thm]{Remark}
\newtheorem{ex}[thm]{Example}

\definecolor{light-gray}{gray}{0.8}
\definecolor{v}{rgb}{0.4,0,0.9}
\definecolor{e}{rgb}{0,1,0.2}
\definecolor{r}{rgb}{1,1,1}

\definecolor{lv}{RGB}{200,180,250}
\definecolor{le}{RGB}{150,255,160}
\definecolor{dy}{RGB}{200, 150, 0}
\definecolor{de}{RGB}{5, 120, 20}

\usepackage{xspace} 

\newcommand{\Z}{\mathbb{Z}}
\newcommand{\HH}{\mathcal H}

\DeclareMathOperator{\bip}{Bip}

\title[Two-variable hypergraph Tutte polynomial via embedding activities]{The two-variable hypergraph Tutte polynomial via embedding activities}

\author{Lilla T\'othm\'er\'esz}
\address{ELTE Eötvös Loránd University, P\'azm\'any P\'eter s\'et\'any 1/C, Budapest, Hungary}
\email{lilla.tothmeresz@ttk.elte.hu}

\begin{document}

\begin{abstract}
	We prove that the two-variable Tutte polynomial of hypergraphs can be defined via embedding activities. We also prove that embedding activities of hypergraphs yield a Crapo-style decomposition of $\mathbb{Z}^E$, thus generalizing Bernardi's results from graphs to hypergraphs.
	
	We also show that hypergraph embedding activities do not fit into the $\Delta$-activity framework of Courtiel. Based on this observation, we construct a graph with an activity notion that yields a Crapo decomposition, but cannot be obtained as a $\Delta$-activity.
\end{abstract}

\maketitle

\section{Introduction}

The Tutte polynomial is one of the most important and well-studied polynomials associated to graphs, and more generally, to matroids. 
One of the nice properties of the Tutte polynomial is that it has many equivalent definitions. A class of these definitions gives the Tutte polynomial as a generating function of various types of activities of spanning trees. The original definition of Tutte used activities with respect to an arbitrary ordering of the edges \cite{Tutte}. A remarkable fact about that definition is that although the activities of individual spanning trees depend on the chosen edge ordering, the generating function does not. A different definition, due to Bernardi \cite{Bernardi_first}, defines activities via an embedding of the graph into an orientable surface, and in this sense, replaces the edge ordering with a more natural auxiliary structure.

Kálmán \cite{hiperTutte} generalized Tutte's activity definition of $T(x,1)$ and $T(1,y)$ to hypergraphs, and more generally, to polymatroids. He called these polynomials the \emph{interior polynomial}, and the \emph{exterior polynomial}, respectively. Later,
Bernardi, Kálmán and Postnikov \cite{BKP} gave a two-variable polynomial $\mathcal{T}_P(x,y)$ for a polymatroid $P$ using activities (again, with respect to a fixed ordering of the ground set). This polynomial generalizes the Tutte polynomial $T_M(x,y)$ of matroids (subject to a change of variables).
A different two-variable polymatroid Tutte polynomial was defined by Fink and Cameron \cite{Cameron-Fink} using lattice point counts instead of activities. Their polynomial is also a transformation of $\mathcal{T}_P(x,y)$.

It is not a priori obvious if Bernardi's embedding activity definition can be generalized further than the class of graphs. However, in \cite{hyperBernardi}, Kálmán and the author of the current paper showed that embedding activities can be generalized to hypergraphs (that can be thought of as a class of polymatroids). They proved that a formula can be given for $\mathcal{T}_\HH(x,1)$ of a hypergraph using embedding activities.
They also conjectured an analogous formula for $\mathcal{T}_\HH(1,y)$ using embedding activities, however, their method of proof could not handle external embedding activities. As a consequence, it also remained open whether a formula can be given for a two-variable hypergraph Tutte polynomial using embedding activities. We note that Bernardi's arguments for the graph case cannot be easily generalized to hypergraphs.

In this paper, we resolve this question, and give a
Bernardi-style definition for the two-variable hypergraph Tutte polynomial. In the mean time, we confirm that the conjectured Bernardi-style definition for $T_\HH(1,y)$ from \cite{hyperBernardi} indeed gives the exterior polynomial. 

We note that \cite{hyperBernardi} also proposed an alternative way to define embedding activities. We conjecture that these alternative embedding activities also yield the polynomial $\mathcal{T}_\HH(x,y)$, but this remains a conjecture (see Section \ref{sec:open}).

In \cite{Courtiel}, Courtiel defined the notion of $\Delta$-activities of graphs (and matroids), which generalizes several notions of activities, including activities with respect to a fixed edge ordering, and embedding activities. He conjectured that if an activity notion (for a graph) yields a Crapo-decomposition, then it fits into the framework of $\Delta$-activities. 

We point out that $\Delta$-activities can be naturally generalized to polymatroids, and they yield a Crapo-decomposition also in this case. However, we show that for hypergraphs, embedding activities cannot always be obtained as $\Delta$-activities, hence the analogue of Courtiel's conjecture is not true for polymatroids, Moreover, based on this observation, we can also construct a graph with an activity notion that yields a Crapo decomposition, but cannot be obtained as a $\Delta$-activity, thus refuting Courtiel's original conjecture.

Let us briefly describe our method. First, we need to recall some more results from \cite{BKP}.
Bernardi, Kálmán and Postnikov \cite{BKP} proved the well-definedness of $\mathcal{T}_P(x,y)$ by first generalizing the corank-nullity definition of the Tutte polynomial to polymatroids. The well-definedness of this corank-nullity polynomial is immediate. Then, they established a Crapo-type decomposition of $\mathbb{Z}^E$ for activities with respect to a fixed edge ordering. Finally, they used their Crapo-decomposition to show that for any fixed ordering, the polynomial defined via activities has the same relationship to the corank-nullity polynomial, hence $\mathcal{T}_P(x,y)$ does not depend on the ordering of the ground set.

In this paper we prove that a Crapo-type decomposition of $\mathbb{Z}^E$ exists with respect to embedding activities of a hypergraph (see Theorem \ref{thm:Crapo_intervals_partition}). Once we prove Theorem \ref{thm:Crapo_intervals_partition}, we can copy \cite{BKP} to prove that the polynomial defined via embedding activities has the same relationship to the corank-nullity polynomial as $\mathcal{T}_P(x,y)$. Hence for any embedding, the polynomial defined via embedding activities agrees with $\mathcal{T}_P(x,y)$.

Outline of the paper: In Section \ref{sec:prelim}, we introduce the necessary background on activities and hypergraphs. We state our main results in Section \ref{sec:corank-nullity}, where we also show how the well-definedness of the embedding activity definition of the two-variable hypergraph Tutte polynomial follows from a Crapo decomposition of $\mathbb{Z}^E$ for embedding activities. Section \ref{sec:existence_of_Crapo} is dedicated to the main technical result: the existence of the Crapo-decomposition for embedding activities of hypergraphs. In Section \ref{sec:open}, we recall an alternative definition for embedding activities suggested in \cite{hyperBernardi}, and pose it as an open problem if those activities also yield $\mathcal{T}_\HH$. 
Finally, in Section \ref{s:Delta_activities}, we investigate Courtiel's $\Delta$-activities. 

\subsection*{Acknowledgement}
I am grateful to Tamás Kálmán for helpful discussions.

This work was supported by the National Research, Development and Innovation Office of Hungary -- NKFIH, grant no.\ 132488, by the János Bolyai Research Scholarship of the Hungarian Academy of Sciences, and by the ÚNKP-22-5, and ÚNKP-23-5 New National Excellence Program of the Ministry for Innovation and Technology, Hungary. This work was also partially supported by the Counting in Sparse Graphs Lendület Research Group of the Alfr\'ed Rényi Institute of Mathematics.

\section{Preliminaries}\label{sec:prelim}

In this section we introduce the necessary background on graphs, hypergraphs and activities.

\subsection{Graphs and their activities}

We briefly recall the definitions of the Tutte polynomial of a graph by Tutte and by Bernardi.

Let $G$ be a graph and $T$ a spanning tree of $G$.

If $e\notin T$, then $T\cup e$ has a unique cycle, which is called the \emph{fundamental cycle} of $e$ in $T$, and is denoted by $C(T,e)$.

If $e\in T$, then $T- e$ has two connected components, and the edges of $G$ connecting the two components form a cut, which is called the \emph{fundamental cut} of $e$ in $T$, and is denoted by $C^*(T,e)$. 

\subsubsection{Activities with respect to an edge ordering}

\begin{defn}[activity (with respect to a fixed edge ordering)]\label{def:graph_activity_fixed_order}
	Let $G$ be a graph, and $<$ an arbitrary ordering of the edges of $G$.
	Let $T$ be a spanning tree of $G$. 
	
	An edge $e\in T$ is internally active if there is no edge $f < e$ such that $T-e\cup f$ is also a spanning tree. $i(T)$ denotes the number of internally active edges for $T$.
	
	An edge $e\notin T$ is externally active if there is no edge $f < e$ such that $T\cup e- f$ is also a spanning tree. $e(T)$ denotes the number of externally active edges for $T$.
\end{defn}

\begin{defn}[Tutte polynomial \cite{Tutte}]
	Let $G$ be a graph.
	$$T_G(x,y)=\sum_{T \text{ spanning tree}} x^{i(T)}y^{e(T)}.$$
\end{defn}
This is well-defined, that is, $T_G(x,y)$ does not depend on the ordering $<$ used to define the activities \cite{Tutte}. 

\subsubsection{Embedding activities}
Let $G$ be a graph. A \emph{ribbon structure} of $G$ is a family of cyclic permutations: for each vertex $x$ of $G$, a cyclic permutation of the edges incident to $x$ is given.
For an edge $xy$ of $G$, we use the following notations: 
\begin{itemize}
	\item $yx^+_G$: the edge following $yx$ at $y$ 
	\item $xy^+_G$: the edge following $xy$ at $x$.
\end{itemize}
If $G$ is clear from the context, we omit the subscript.

In addition to a ribbon structure, we also need to fix a \emph{basis} $(b_0,b_0b_1)$, where $b_0$ is an arbitrary node of the graph and $b_0b_1$ is an arbitrary edge incident to $b_0$. 

Suppose that a ribbon structure and a basis are fixed. Then any spanning tree $T$ of $G$ gives a ``walk'' in the graph. This was defined by Bernardi \cite{Bernardi_first}, and following him we call it the \emph{tour} of $T$.

\begin{defn}[Tour of a spanning tree] \label{def:tour_of_a_tree}
	Let $G$ be a ribbon graph with a basis $(b_0,b_0b_1)$, and let $T$ be a spanning tree of $G$.
	The \emph{tour} of $T$ is a sequence of node-edge pairs, starting with $(b_0, b_0b_1)$. If the current node-edge pair is $(x,xy)$ and $xy\notin T$, then the current node-edge pair of the next step is $(x,xy^+)$. If the current node-edge pair is $(x,xy)$ and $xy\in T$, then the current node-edge pair of the next step is $(y,yx^+)$.	The tour stops right before when $(b_0,b_0b_1)$ would once again become the current node-edge pair. 
\end{defn}

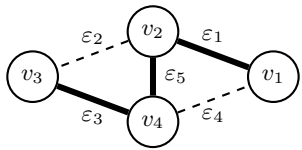
\begin{figure} 
	\begin{center}
		\begin{tikzpicture}[-,>=stealth',auto,scale=0.4,
			thick]
			\tikzstyle{o}=[circle,draw]
			\node[o] (1) at (8, 0) {{\small $v_1$}};
			\node[o] (2) at (4, 1.5) {{\small $v_2$}};
			\node[o] (3) at (0, 0) {{\small $v_3$}};
			\node[o] (4) at (4, -1.5) {{\small $v_4$}};
			\path[-,every node/.style={font=\sffamily\small}, line width=0.8mm]
			(1) edge node [above] {$\varepsilon_1$} (2)
			(3) edge node [below] {$\varepsilon_3$} (4)
			(2) edge node {$\varepsilon_5$} (4);
			\path[-,every node/.style={font=\sffamily\small},dashed]
			(4) edge node [below] {$\varepsilon_4$} (1)
			(2) edge node [above] {$\varepsilon_2$} (3);
		\end{tikzpicture}
	\end{center}
	\caption{The tour of a spanning tree. Let the ribbon structure be the one induced by the positive orientation of the plane. The edges of the tree are drawn by thick lines, the non-edges by dashed lines. With $b_0=v_1, b_1=v_2$, we get the tour $(v_1,\varepsilon_1)$, $(v_2,\varepsilon_2)$, $(v_2,\varepsilon_5)$, $(v_4,\varepsilon_3)$, $(v_3,\varepsilon_2)$, $(v_3,\varepsilon_3)$, $(v_4,\varepsilon_4)$, $(v_4,\varepsilon_5)$, $(v_2,\varepsilon_1)$, $(v_1,\varepsilon_4)$. 
	}
	\label{fig:tour_of_a_tree}
\end{figure} 

See Figure \ref{fig:tour_of_a_tree} for an example. 
Bernardi proved {\cite[Lemma 4.2 (Lemma 5 in the arxiv version)]{Bernardi_first}} that in the tour of a spanning tree $T$, each edge $xy$ of $G$ becomes current edge twice, in one case with $x$ as current vertex, and in the other case with $y$ as current vertex.
This naturally orders the edges of the graph in the order that they first become current.

Bernardi defined the internal and external embedding activities of a tree $T$ as the internal and external activities with respect to the order induced by the tour of $T$. He then showed that the two-variable generating function of internal and external activities also gives the Tutte polynomial (in particular, the generating function does not depend on the chosen ribbon structure and basis).

Embedding activities were generalized for hypergraphs in \cite{hyperBernardi}, where also their connections with geometry were explored. (They are related to dissections of root polytopes of bipartite graphs.) We give these definitions in the next section. 
In \cite{hyperBernardi}, a formula for $\mathcal{T}_\HH(x,1)$ was given using embeding activities. In this paper, we show that $\mathcal{T}_\HH(x,y)$ can also be defined via embedding activities.

\subsection{Polymatroids, hypergraphs and their activities}

\subsubsection{Polymatroids}
A polymatroids were introduced by Edmonds \cite{Edmonds_polymatroid}, as a generalization of matroids. Here, by polymatroids, we will mean what is commonly called an integer polymatroid.

\begin{defn}\label{def:polymatroid}
	Let $P\subset \mathbb{Z}^E$.
	For any $b, b' \in P$, if for some $e\in E$ we have $b(e) < b'(e)$, then
	there exist some $f\in E$ with $b(f) > b'(f)$ such that $b+\mathbf{1}_{e}-\mathbf{1}_{f}$ and $b'-\mathbf{1}_{e}+\mathbf{1}_{f}$ are both in $P$.
\end{defn}

We note that Edmonds also requires that $P\subseteq \mathbb{Z}^V_{\geq 0}$. We follow \cite{BKP} in omitting this requirement, as the definition of the polymatroid Tutte polynomial is translation invariant.

An important example for a polymatroid is a matroid: If one takes the characteristic vectors of bases, then the symmetric exchange theorem implies that we get a polymatroid in the sense of Definition \ref{def:polymatroid}. Another important example for polymatroids is hypergraphs, that we elaborate in section \ref{ss:hypergraph}.

\subsubsection{Polymatroid activities with respect to a fixed ordering}
Let us recall how \cite{hiperTutte} and \cite{BKP} associates activities to bases of a polymatroid. 

\begin{defn}
	\label{def:activity_w_r_to_order}
	Let $P$ be a polymatroid on the ground set $E$ with an order $<$ on $E$. 
	
	An element $e\in E$ is \emph{internally active} for the basis $b$, with respect to $<$, if $b-\mathbf{1}_e +\mathbf{1}_f$ is not a basis for any $f<e$. 
	
	An element $e\in E$ is \emph{externally active} for the basis $b$, with respect to $<$, if $b+\mathbf{1}_e -\mathbf{1}_f$ is not a basis for any $f<e$.
\end{defn}

\begin{remark}\label{rem:graph_polymatroid_activity_diff}
We note activities for graphs (Definition \ref{def:graph_activity_fixed_order}) are not a direct special case of this definition, as Definition \ref{def:graph_activity_fixed_order}  required that $e\in T$ for an internally active edge, and $e\notin T$ for an externally active edge. As polmatroids can take values other than 0 and 1, those conditions cannot really be generalized to the polymatroid case, however, by Definition \ref{def:activity_w_r_to_order} an edge $e\notin T$ is always internally active, and an edge $e\in T$ is always externally active. 
Hence compared to the graph activity definition, Definition \ref{def:activity_w_r_to_order} shifts the internal activity of each spanning tree by $|V|-1$ and the external activity of each spanning tree by $|E|-|V|+1$, modifying the generating function in a trivial way. 
\end{remark}

Let $Int_<(b)$ denote the set of internally active elements of $b$ with respect to the order $<$, and let $Ext_<(b)$ denote the set of externally active elements of $b$ with respect to $<$.

Let $i_<(b)=|Int_<(b)|$ denote the number of internally active elements in $b$ and call this value the \emph{internal activity} of $b$. 
Let $e_<(b)=|Ext_<(b)|$ denote the number of externally active elements in $b$ and call this value the \emph{external activity} of $b$.

Moreover, let $oi_<(b)=|Int_<(b)-Ext_<(b)|$, $oe_<(b)=|Ext_<(b)-Int_<(b)|$ and $ie_<(b)=|Int_<(b)\cap Ext_<(b)|$ be respectively the number of only internally active, only externally active, and both internally and externally active elements.

The polymatroid Tutte polynomial of Bernardi, Kálmán and Postnikov is defined the following way:
\begin{defn}\label{def:Tutte_poly_activity_fixed_ordering}
	$$\mathcal{T}_P(x,y)=\sum_{b\in P} x^{oi_<(b)}y^{oe_<(b)}(x+y-1)^{ie_<(b)}.$$
\end{defn}
They show that this polynomial is well-defined, that is, it does not depend on the chosen ordering of $E$. Moreover, it has many nice properties. 
The polymatroid Tutte polynomial does not agree with the usual graphic Tutte polynomial for graphic matroids, but they have the following simple relationship (see \cite[Theorem 5.2]{BKP}):
$T_G(x,y)=x^{|E|-|V|+1}y^{|V|-1} \mathcal{T}_G(\frac{x+y-1}{y},\frac{x+y-1}{y})$.

Two interesting specializations of $\mathcal{T}_P(x,y)$ are
$$\mathcal{T}_P(x,1)=\sum_{b\in P} x^{i_<(b)} \quad\text{ and }\quad \mathcal{T}_P(1,y)=\sum_{b\in P} y^{e_<(b)},$$ that we call the \emph{interior} and \emph{exterior} polynomial, respectively. We note that in \cite{hiperTutte}, $x^{|E|}\mathcal{T}_P(\frac{1}{x},1)$ is called the interior polynomial, and 
$y^{|E|}\mathcal{T}_P(1,\frac{1}{y})$ is called the exterior polynomial.

\subsubsection{Hypergraphs}\label{ss:hypergraph}
A \emph{hypergraph} is an ordered pair $\HH=(V,E)$, where $V$ is a finite set and $E$ is a finite multiset of subsets of $V$. We refer to elements of $V$ as \emph{vertices} and to elements of $E$ as \emph{hyperedges}. 
It is convenient to represent a hypergraph by a bipartite graph.
For a hypergraph $\mathcal{H}=(V,E)$, let the \emph{underlying bipartite graph} $\bip\HH$ be the bipartite graph with vertex classes $V$ and $E$, where $v\in V$ is connected to $e\in E$ if $v\in e$ in $\HH$. 
We will mostly think of hypergraphs as bipartite graphs.
When talking about $\bip\HH$, we will call the elements of $V\cup E$ \emph{nodes}. Specifically, the elements of $V$ will be called \emph{violet} nodes and elements of $E$ \emph{emerald} nodes. We use the words hyperedge and emerald node interchangeably. We use greek letters to denote the edges of $\bip\HH$.

Throughout the paper, we assume that $\HH$ is \emph{connected}, by which we mean that $\bip \HH$ is connected.

Spanning trees of graphs are the bases of a matroid: the graphic matroid. For hypergraphs, Kálmán \cite{hiperTutte} defined hypertrees, that generalize (characteristic vectors of) spanning trees from graphs to hypergraphs. 

\begin{defn}
	\label{def:hypertree}
	Let $\mathcal{H}=(V,E)$ be a hypergraph with underlying bipartite graph $\bip\HH$. We say that the vector $h\in \Z^E$ is a  \emph{hypertree} if there exists a spanning tree $T$ of $\bip\HH$ that has degree $d_T(e)=h(e)+1$ at each node $e\in E$. In this case we say that $T$ represents $h$.
\end{defn}
Most importantly for us, hypertrees form the bases of a polymatroid \cite{hiperTutte}.

\begin{prop}\cite{hiperTutte}\label{prop:exchange_property_for_hypertrees}
	Let $h$ and $h'$ be hypertrees and $e$ be a hyperedge with $h(e) < h'(e)$. Then
	there exist a hyperedge $f$ with $h(f) > h'(f)$ such that $h+\mathbf{1}_{e}-\mathbf{1}_{f}$ and $h'-\mathbf{1}_{e}+\mathbf{1}_{f}$ are both hypertrees.
\end{prop}

We denote the set of all hypertrees of $\HH$ by $H(\HH)$. If $\HH$ is clear from the context, we simply write $H$.

Note that a hypertree might have many different representing spanning trees.

It is an easy exercise to check that if $G$ is a connected graph, then the hypertrees of $\bip G$ are exactly the characteristic vectors of spanning trees, hence hypertrees indeed generalize spanning trees.

\subsubsection{Embedding activities for hypergraphs}

To define embedding activities for hypergraphs, we assume that $\bip\HH$ has a ribbon graph structure, and a basis $(b_0, b_0b_1)$ is fixed (where $b_0$ is a node of $\bip \HH$ (either from $V$ or from $E$) and $b_0b_1$ is an edge of $\bip \HH$ incident to $b_0$).

To generalize Bernardi's embedding activities, one needs to generalize the tour of a spanning tree to hypertrees. To do this, we will make use of a nice representation of hypertrees.

Recall that hypertrees are exactly those vectors $h\in\mathbb{Z}^E$ such that $h+\mathbf{1}_E$ is the degree sequence of some spanning tree of $G$. It is very useful to have a system of representing spanning trees at hand. Fortunately, ribbon structures provide us such trees.

\begin{defn}[Jaeger tree \cite{hyperBernardi}]
	In a bipartite graph with a ribbon structure and basis, we call a spanning tree $T$ a \emph{Jaeger tree} if for each edge $ev \notin T$ with $e\in E$ and $v\in V$, the tour of $T$ has $(e,ev)$ as a current node-edge pair before $(v,ev)$. In other words, in the tour of a Jaeger tree, each non-edge is first seen at its emerald endpoint.
\end{defn}

\begin{figure}[t]
	\begin{tikzpicture}[scale=.23]
		
		\begin{scope}[shift={(-14,0)}]
			\draw [ultra thick] (18, 8.5) -- (14,6.5);
			\draw [ultra thick] (14, 6.5) -- (14,2);
			\draw [ultra thick] (18, 8.5) -- (22,6.5);
			\draw [dashed, thick] (18, 8.5) -- (18,4);
			\draw [ultra thick] (14, 2) -- (18,4);
			\draw [ultra thick] (14, 2) -- (18,0);
			\draw [dashed, thick] (18, 0) -- (22,2);
			\draw [dashed, thick] (22, 2) -- (18,4);
			\draw [ultra thick] (22, 2) -- (22,6.5);
			\draw [fill=e,e] (18, 0) circle [radius=0.6];
			\draw [fill=e,e] (18, 4) circle [radius=0.6];
			\draw [fill=e,e] (14, 6.5) circle [radius=0.6];
			\draw [fill=e,e] (22, 6.5) circle [radius=0.6];
			\draw [fill=v,v] (14, 2) circle [radius=0.6];
			\draw [fill=v,v] (22, 2) circle [radius=0.6];
			\draw [fill=v,v] (18, 8.5) circle [radius=0.6];
			
			\draw [fill] (12.5,1.5) circle [radius=.3];
			\draw [->] (12.5,1.5) -- (13.6, 0.9);
			
			\node at (14, 6.5) {{\tiny{$1$}}};
			\node at (22, 6.5) {{\tiny{$1$}}};
			\node at (18, 4) {{\tiny{$0$}}};
			\node at (18, 0) {{\tiny{$0$}}};
			
			\node at (18, -2.2) {{\small{$y^2(x+y-1)^2$}}};
		\end{scope}
		
		\begin{scope}[shift={(0,0)}]
			\draw [ultra thick] (18, 8.5) -- (14,6.5);
			\draw [dashed, thick] (14, 6.5) -- (14,2);
			\draw [ultra thick] (18, 8.5) -- (22,6.5);
			\draw [ultra thick] (18, 8.5) -- (18,4);
			\draw [ultra thick] (14, 2) -- (18,4);
			\draw [ultra thick] (14, 2) -- (18,0);
			\draw [dashed, thick] (18, 0) -- (22,2);
			\draw [dashed, thick] (22, 2) -- (18,4);
			\draw [ultra thick] (22, 2) -- (22,6.5);
			\draw [fill=e,e] (18, 0) circle [radius=0.6];
			\draw [fill=e,e] (18, 4) circle [radius=0.6];
			\draw [fill=e,e] (14, 6.5) circle [radius=0.6];
			\draw [fill=e,e] (22, 6.5) circle [radius=0.6];
			\draw [fill=v,v] (14, 2) circle [radius=0.6];
			\draw [fill=v,v] (22, 2) circle [radius=0.6];
			\draw [fill=v,v] (18, 8.5) circle [radius=0.6];
			
			\node at (14, 6.5) {{\tiny{$0$}}};
			\node at (22, 6.5) {{\tiny{$1$}}};
			\node at (18, 4) {{\tiny{$1$}}};
			\node at (18, 0) {{\tiny{$0$}}};
			
			\node at (18, -2.2) {{\small{$xy^2(x+y-1)$}}};
		\end{scope}
		
		\begin{scope}[shift={(14, 0)}]
			\draw [ultra thick] (18, 8.5) -- (14,6.5);
			\draw [ultra thick] (14, 6.5) -- (14,2);
			\draw [dashed, thick] (18, 8.5) -- (22,6.5);
			\draw [dashed, thick] (18, 8.5) -- (18,4);
			\draw [ultra thick] (14, 2) -- (18,4);
			\draw [ultra thick] (14, 2) -- (18,0);
			\draw [dashed, thick] (18, 0) -- (22,2);
			\draw [ultra thick] (22, 2) -- (18,4);
			\draw [ultra thick] (22, 2) -- (22,6.5);
			\draw [fill=e,e] (18, 0) circle [radius=0.6];
			\draw [fill=e,e] (18, 4) circle [radius=0.6];
			\draw [fill=e,e] (14, 6.5) circle [radius=0.6];
			\draw [fill=e,e] (22, 6.5) circle [radius=0.6];
			\draw [fill=v,v] (14, 2) circle [radius=0.6];
			\draw [fill=v,v] (22, 2) circle [radius=0.6];
			\draw [fill=v,v] (18, 8.5) circle [radius=0.6];
			
			\node at (14, 6.5) {{\tiny{$1$}}};
			\node at (22, 6.5) {{\tiny{$0$}}};
			\node at (18, 4) {{\tiny{$1$}}};
			\node at (18, 0) {{\tiny{$0$}}};
			
			\node at (18, -2.2) {{\small{$xy^2(x+y-1)$}}};
		\end{scope}
		
		\begin{scope}[shift={(28, 0)}]
			\draw [ultra thick] (18, 8.5) -- (14,6.5);
			\draw [dashed, thick] (14, 6.5) -- (14,2);
			\draw [dashed, thick] (18, 8.5) -- (22,6.5);
			\draw [ultra thick] (18, 8.5) -- (18,4);
			\draw [ultra thick] (14, 2) -- (18,4);
			\draw [ultra thick] (14, 2) -- (18,0);
			\draw [dashed, thick] (18, 0) -- (22,2);
			\draw [ultra thick] (22, 2) -- (18,4);
			\draw [ultra thick] (22, 2) -- (22,6.5);
			\draw [fill=e,e] (18, 0) circle [radius=0.6];
			\draw [fill=e,e] (18, 4) circle [radius=0.6];
			\draw [fill=e,e] (14, 6.5) circle [radius=0.6];
			\draw [fill=e,e] (22, 6.5) circle [radius=0.6];
			\draw [fill=v,v] (14, 2) circle [radius=0.6];
			\draw [fill=v,v] (22, 2) circle [radius=0.6];
			\draw [fill=v,v] (18, 8.5) circle [radius=0.6];
			
			\node at (14, 6.5) {{\tiny{$0$}}};
			\node at (22, 6.5) {{\tiny{$0$}}};
			\node at (18, 4) {{\tiny{$2$}}};
			\node at (18, 0) {{\tiny{$0$}}};
			
			\node at (18, -2.2) {{\small{$x^2y(x+y-1)$}}};
		\end{scope}
		
		\begin{scope}[shift={(-14, -14)}]
			\draw [ultra thick] (18, 8.5) -- (14,6.5);
			\draw [ultra thick] (14, 6.5) -- (14,2);
			\draw [dashed, thick] (18, 8.5) -- (22,6.5);
			\draw [dashed, thick] (18, 8.5) -- (18,4);
			\draw [dashed, thick] (14, 2) -- (18,4);
			\draw [ultra thick] (14, 2) -- (18,0);
			\draw [ultra thick] (18, 0) -- (22,2);
			\draw [ultra thick] (22, 2) -- (18,4);
			\draw [ultra thick] (22, 2) -- (22,6.5);
			\draw [fill=e,e] (18, 0) circle [radius=0.6];
			\draw [fill=e,e] (18, 4) circle [radius=0.6];
			\draw [fill=e,e] (14, 6.5) circle [radius=0.6];
			\draw [fill=e,e] (22, 6.5) circle [radius=0.6];
			\draw [fill=v,v] (14, 2) circle [radius=0.6];
			\draw [fill=v,v] (22, 2) circle [radius=0.6];
			\draw [fill=v,v] (18, 8.5) circle [radius=0.6];
			
			\node at (14, 6.5) {{\tiny{$1$}}};
			\node at (22, 6.5) {{\tiny{$0$}}};
			\node at (18, 4) {{\tiny{$0$}}};
			\node at (18, 0) {{\tiny{$1$}}};
			
			\node at (18, -2.2) {{\small{$x^2y(x+y-1)$}}};
		\end{scope}
		
		\begin{scope}[shift={(0, -14)}]
			\draw [ultra thick] (18, 8.5) -- (14,6.5);
			\draw [dashed, thick] (14, 6.5) -- (14,2);
			\draw [dashed, thick] (18, 8.5) -- (22,6.5);
			\draw [ultra thick] (18, 8.5) -- (18,4);
			\draw [dashed, thick] (14, 2) -- (18,4);
			\draw [ultra thick] (14, 2) -- (18,0);
			\draw [ultra thick] (18, 0) -- (22,2);
			\draw [ultra thick] (22, 2) -- (18,4);
			\draw [ultra thick] (22, 2) -- (22,6.5);
			\draw [fill=e,e] (18, 0) circle [radius=0.6];
			\draw [fill=e,e] (18, 4) circle [radius=0.6];
			\draw [fill=e,e] (14, 6.5) circle [radius=0.6];
			\draw [fill=e,e] (22, 6.5) circle [radius=0.6];
			\draw [fill=v,v] (14, 2) circle [radius=0.6];
			\draw [fill=v,v] (22, 2) circle [radius=0.6];
			\draw [fill=v,v] (18, 8.5) circle [radius=0.6];
			
			\node at (14, 6.5) {{\tiny{$0$}}};
			\node at (22, 6.5) {{\tiny{$0$}}};
			\node at (18, 4) {{\tiny{$1$}}};
			\node at (18, 0) {{\tiny{$1$}}};
			
			\node at (18, -2.2) {{\small{$x^2(x+y-1)$}}};
		\end{scope}
		
		\begin{scope}[shift={(14, -14)}]
			\draw [ultra thick] (18, 8.5) -- (14,6.5);
			\draw [dashed, thick] (14, 6.5) -- (14,2);
			\draw [ultra thick] (18, 8.5) -- (22,6.5);
			\draw [ultra thick] (18, 8.5) -- (18,4);
			\draw [dashed, thick] (14, 2) -- (18,4);
			\draw [ultra thick] (14, 2) -- (18,0);
			\draw [ultra thick] (18, 0) -- (22,2);
			\draw [dashed, thick] (22, 2) -- (18,4);
			\draw [ultra thick] (22, 2) -- (22,6.5);
			\draw [fill=e,e] (18, 0) circle [radius=0.6];
			\draw [fill=e,e] (18, 4) circle [radius=0.6];
			\draw [fill=e,e] (14, 6.5) circle [radius=0.6];
			\draw [fill=e,e] (22, 6.5) circle [radius=0.6];
			\draw [fill=v,v] (14, 2) circle [radius=0.6];
			\draw [fill=v,v] (22, 2) circle [radius=0.6];
			\draw [fill=v,v] (18, 8.5) circle [radius=0.6];
			
			\node at (14, 6.5) {{\tiny{$0$}}};
			\node at (22, 6.5) {{\tiny{$1$}}};
			\node at (18, 4) {{\tiny{$0$}}};
			\node at (18, 0) {{\tiny{$1$}}};
			
			\node at (18, -2.2) {{\small{$x^2(x+y-1)^2$}}};
		\end{scope}
		\begin{scope}[shift={(28, -14)}]
			\draw [thick] (18, 8.5) -- (14,6.5);
			\draw [thick] (14, 6.5) -- (14,2);
			\draw [thick] (18, 8.5) -- (22,6.5);
			\draw [thick] (18, 8.5) -- (18,4);
			\draw [thick] (14, 2) -- (18,4);
			\draw [thick] (14, 2) -- (18,0);
			\draw [thick] (18, 0) -- (22,2);
			\draw [thick] (22, 2) -- (18,4);
			\draw [thick] (22, 2) -- (22,6.5);
			\draw [fill, le] (18, 0) circle [radius=0.7];
			\draw [fill, le] (18, 4) circle [radius=0.7];
			\draw [fill, le] (14, 6.5) circle [radius=0.7];
			\draw [fill, le] (22, 6.5) circle [radius=0.7];
			\draw [fill, lv] (14, 2) circle [radius=0.7];
			\draw [fill, lv] (22, 2) circle [radius=0.7];
			\draw [fill, lv] (18, 8.5) circle [radius=0.7];
			
			\node at (14, 6.5) {{\tiny{$e_3$}}};
			\node at (22, 6.5) {{\tiny{$e_2$}}};
			\node at (18, 4) {{\tiny{$e_1$}}};
			\node at (18, 0) {{\tiny{$e_0$}}};
			
			\node at (14, 2) {{\tiny{$v_0$}}};
			\node at (22, 2) {{\tiny{$v_1$}}};
			\node at (18, 8.5) {{\tiny{$v_2$}}};
			
	    \end{scope}
	\end{tikzpicture}
	\caption{A bipartite graph corresponding to a hypergraph with 3 vertices (violet nodes) and 4 hyperedges (emerald nodes). Take the ribbon structure induced by the positive orientation of the plane, and basis $(v_0,v_0e_0)$ (see the last panel for the notations). The first 7 panels show the 7 hypertrees of the hypergraph, together with their representing Jaeger trees, and the corresponding terms of $\mathcal{T}^{em}_\HH$.}\label{fig:embedding_activity}
\end{figure}
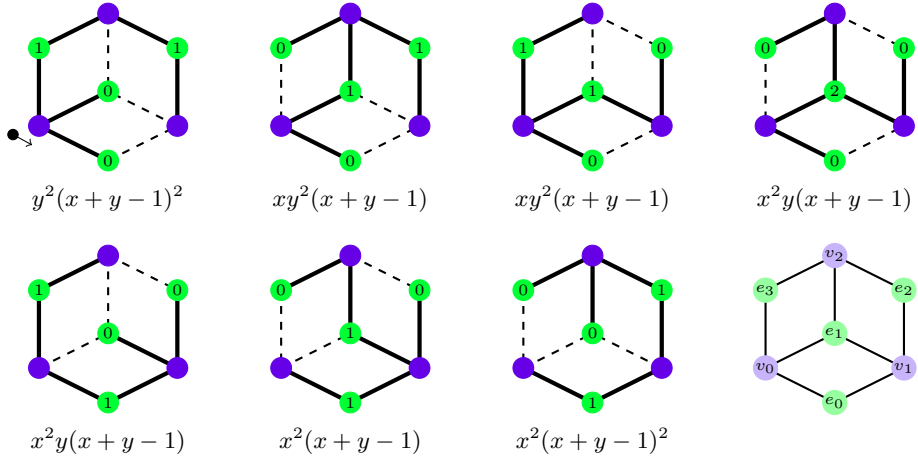

\begin{ex}
	See Figure \ref{fig:embedding_activity} for examples of Jaeger trees. The last panel shows the labels of the nodes. In particular, the tour of the tree on the first panel is $(v_0,v_0e_0)$, $(e_0,e_0v_1)$, $(e_0,e_0v_0)$, $(v_0,v_0e_1)$, $(e_1,e_1v_1)$,  $(e_1,e_1v_2)$, $(e_1,e_1v_0)$,  $(v_0,v_0e_3)$, $(e_3,e_3v_2)$, $(v_2,v_2e_1)$, $(v_2,v_2e_2)$, $(e_2,e_2v_1)$, $(v_1,v_1e_1)$, $(v_1,v_1e_0)$, $(v_1,v_1e_2)$, $(e_2,e_2v_2)$, $(v_2,v_2e_3)$, $(e_3,e_3v_0)$. This is a Jaeger tree since $(e_0,e_0v_1)$ precedes $(v_1,v_1e_0)$, $(e_1,e_1v_1)$ precedes $(v_1,v_1e_1)$ and $(e_1,e_1v_2)$ precedes $(v_2,v_2e_1)$.
\end{ex}
\begin{thm}\cite{hyperBernardi}
	Let $\HH$ be a connected hypergraph, and fix an arbitrary ribbon structure and basis.
	Then for each hypertree $h$, there is exactly one Jaeger tree representing $h$.
\end{thm}

\begin{remark}\label{rem:emerald_vs_violet_Jaeger_tree}
	Of course one could also define Jaeger trees with the rule that each non-tree edge is first seen at its violet endpoint. In \cite{hyperBernardi} both types of trees are investigated, and are called emerald, and violet Jaeger trees, respectively. In this paper, we only consider emerald Jaeger trees, and hence we drop the word ``emerald''. 
	
We note that Jaeger trees also have a nice geometric interpretation as a dissection of the root polytope of $\bip\HH$ (see \cite[Section 6]{hyperBernardi}). While this interpretation is often very useful, it will not play a role in this paper.
\end{remark}

Since given a ribbon structure and basis, now we have a well-defined set of representing spanning trees for the hypertrees, we can associate to hypertrees the tours of their Jaeger trees. Then, we can use the tour to define an ordering of the hyperedges, as we now explain.

\begin{defn}[ordering corresponding to a hypertree, $<_h$]
	Let $\HH=(V,E)$ be a hypergraph and fix a ribbon structure on $\bip\HH$ and a basis $(b_0,b_0b_1)$. Let $h$ be a hypertree. Then $<_h$ is a complete ordering of $E$ defined the following way. 
	
	Let $T$ be the unique Jaeger tree realizing $h$. 
	For $e,f\in E$, we define $e <_h f$, if in the tour of $T$, the first time that $e$ is reached is before the first time that $f$ is reached. We denote $e \leq_h f$ if either $e=f$ or $e <_h f$.
\end{defn}
\begin{ex}
	The first panel of Figure \ref{fig:embedding_activity} shows the hyppertree $h(e_0)=0$, $h(e_1)=0$, $h(e_2)=1$, $h(e_3)=1$. The ordering $<_h$ corresponding to this hypertree is $e_0 < e_1 < e_3 < e_2$. 
\end{ex}

\begin{defn}[Embedding activities]
	\label{def:embedding_activity}
	Let $\HH=(V,E)$ be a hypergraph and fix a ribbon structure on $\bip\HH$ and a basis $(b_0,b_0b_1)$.
	
	A hyperedge $e\in E$ is \emph{internally embedding active} for the hypertree $h$ (with respect to the ribbon structure and basis), if $h-\mathbf{1}_e +\mathbf{1}_f$ is not a hypertree for any $f <_h e$. 
	
	A hyperedge $e\in E$ is \emph{externally embedding active} for the hypertree $h$ (with respect to the ribbon structure and basis), if $h+\mathbf{1}_e -\mathbf{1}_f$ is not a hypertree for any $f <_h e$.
\end{defn}

We denote by $Int_{em}(h)$ the set of internally embedding active elements of $h$, and by $Ext_{em}(h)$ the set of externally embedding active elements of $h$.

We denote $oi_{em}(h)=|Int_{em}(h)-Ext_{em}(h)|$, $oe_{em}(h)=|Ext_{em}(h)-Int_{em}(h)|$, $ie_{em}(h)=|Int_{em}(h)\cap Ext_{em}(h)|$.

\begin{ex}
	Take the graph $\bip \HH$ of Figure \ref{fig:embedding_activity} with the indicated ribbon structure and basis. (For the node names, see the last panel.) 
	The $6^{th}$ panel shows the Jaeger tree of the hypertree $h$ with $h(e_0)=h(e_1)=1$, $h(e_2)=h(e_3)=0$. We have $e_0 <_h e_2 <_h e_1 <_h e_3$. For the hypertree $h$, hyperedge $e_0$ is both internally and externally active (since it is the smallest element with respect to $<_h$). Hyperedge $e_2$ is internally active, since a hypertree is nonnegative, however, it is externally passive since $h+\mathbf{1}_{e_2}-\mathbf{1}_{e_0}$ is a hypertree. Hyperedge $e_1$ is both internally and externally passive, since $h+\mathbf{1}_{e_1}-\mathbf{1}_{e_0}$ and $h-\mathbf{1}_{e_1}+\mathbf{1}_{e_2}$ are both hypertrees. Finally, $e_3$ is internally active since $h(e_3)=0$, but externally passive since $h+\mathbf{1}_{e_3}-\mathbf{1}_{e_0}$ is a hypertree.
\end{ex}

\section{Well-definedness of the two variable embedding Tutte polynomial}\label{sec:corank-nullity}
Let us define the two-variable hypergraph Tutte polynomial via embedding activities, analogously to the formula using activities with respect to an ordering of the hyperegdes.
\begin{defn}\label{def:hypergraph_embedding_Tutte_poly}
	$$\mathcal{T}^{em}_\HH(x,y)=\sum_{h\in H} x^{oi_{em}(h)}y^{oe_{em}(h)}(x+y-1)^{ie_{em}(h)}.$$
\end{defn}

\begin{ex}
	For the hypergraph of Figure \ref{fig:embedding_activity}, we have
	$$\mathcal{T}^{em}_\HH(x,y)=
	\begin{array}{ccccc}
	&  &  y^2&  -2y^3&  +y^4\\  
	&  &  -4xy^2&  +4xy^3&  \\  
	&  -3x^2y&  +6x^2y^2&  &  \\  
	-x^3&  +4x^3y&  &  &  \\  
	+x^4&  &  &  &  
   \end{array}.$$
\end{ex}

\begin{thm}\label{thm:embedding_Tutte_well_defined}
	For any hypergraph $\HH$, $\mathcal{T}^{em}_\HH$ is well-defined, that is, it does not depend on the chosen ribbon structure and basis. Moreover, $\mathcal{T}^{em}_\HH=\mathcal{T}_\HH$.
\end{thm}

We prove Theorem \ref{thm:embedding_Tutte_well_defined} later in this Section.

Note the following two specializatons: $$\mathcal{T}^{em}_\HH(x,1)=\sum_{h\in H} x^{oi_{em}(h)}x^{ie_{em}(h)}=x^{|Int_{em}(h)|}$$ and 
$$\mathcal{T}^{em}_\HH(1,y)=\sum_{h\in H} y^{oe_{em}(h)}y^{ie_{em}(h)}=y^{|Ext(h)|}.$$
	
	In \cite{hyperBernardi}, Kálmán and the author proved that $\mathcal{T}^{em}_\HH(x,1)$ is well-defined and agrees with $\mathcal{T}_\HH(x,1)$. Theorem \ref{thm:embedding_Tutte_well_defined} also implies the following, which was stated as a conjecture in \cite{hyperBernardi}.
	
	\begin{coroll}
		$\mathcal{T}^{em}_\HH(1,y)$ is well-defined and it agrees with $\mathcal{T}_\HH(1,y)$.
	\end{coroll}

To prove Theorem \ref{thm:embedding_Tutte_well_defined}, we recall the definition of the corank-nullity polynomial $\tilde{\mathcal{T}}_\HH$ of a hypergraph from \cite{BKP}, and show that $\mathcal{T}^{em}_\HH$ has the same relationship to $\tilde{\mathcal{T}}_\HH$ as 
$\mathcal{T}_\HH$ does, thus proving Theorem \ref{thm:embedding_Tutte_well_defined}.

To establish the relationship of $\mathcal{T}^{em}_\HH$ and $\tilde{\mathcal{T}}_\HH$, we will copy \cite{BKP}: They proved the existence of a Crapo-type decomposition of $\mathbb{Z}^E$ for activities with respect to a hyperedge ordering, then used the Crapo decomposition to relate $\mathcal{T}_\HH$ to $\tilde{\mathcal{T}}_\HH$.
We will apply the same argument, with the only exception that we need to prove the existence of a Crapo decomposition for embedding activities. The proof of the embedding Crapo decomposition is the main result of this paper, and it is deferred to Section \ref{sec:existence_of_Crapo}.
Let us introduce the necessary notions.

\begin{defn}[Embedding Crapo interval]  
	The Crapo interval of a hypertree $h$ with respect to a ribbon structure and basis is 
	\begin{align*}
	C_{em}(h)=\{ y \in \mathbb{Z}^E : (y(e) > h(e) \Rightarrow e \in Ext_{em}(h)), (y(e) < h(e) \Rightarrow e\in Int_{em}(h))\}
	\end{align*}
\end{defn}

For $c\in\mathbb{Z}^E$ and $S\subseteq \mathbb{Z}^E$ let us denote $d_1(S,c)=\min\{\sum_{e\in E}|s(e)-c(e)| \mid s\in S\}$.
This is commonly called the Manhattan distance.

The following is the main technical theorem of the paper.

\begin{thm}\label{thm:Crapo_intervals_partition} 
\begin{enumerate}
	\item\label{st:partition} The Crapo intervals $\{C_{em}(h) \mid h\in H(\HH)\}$ partition $\mathbb{Z}^E$.	
	
	\item\label{st:min} For each $c\in\mathbb{Z}^E$, if $c\in C_{em}(h)$, then $d_1(H,c)=d_1(h,c)$.
\end{enumerate}
\end{thm}

We give the proof in Section \ref{sec:existence_of_Crapo}.

We also need the generalizations of corank and nullity from \cite{BKP}:

\begin{defn}\cite{BKP}
	 Let $S\subseteq \mathbb{Z}^E$ and $c\in\mathbb{Z}^E$. We denote \\
	$d_1^<(S,c)=\min\{\sum_{e\in E} \max\{0, c(e)-s(e)\} \mid s\in S\}$,\\
	$d_1^>(S,c)=\min\{\sum_{e\in E} \max\{0, s(e)-c(e)\} \mid s\in S\}$.
\end{defn}

Here
$d_1^<(H,c)$ is a generalization of corank. Indeed, if we take a graph, that is, the set of hypertrees $H$ is the set of characteristic vectors of spanning trees, then for a $c\in\{0,1\}^E$, $d_1(H,c)$ is exactly the corank of $c$. Similarly, $d_1^>(H,c)$ is the generalization of nullity.

Using these notions, Bernardi, Kálmán and Postnikov define the following corank-nullity polynomial:
$$\tilde{\mathcal{T}}_\HH(u,v)=\sum_{c\in\mathbb{Z}^E} u^{d_1^>(H(\HH),c)}v^{d_1^<(H(\HH),c)}$$

It is clear that $\tilde{\mathcal{T}}_\HH$ is well-defined. Bernardi, Kálmán and Postnikov \cite{BKP} shows the well-definedness of $\mathcal{T}_\HH$ by showing that 
$$ \tilde{\mathcal{T}}_\HH(u,v) = \mathcal{T}_\HH\left(\frac{1}{1-u}, \frac{1}{1-v} \right)
$$
which is understood as an identity of formal power series, where $\frac{1}{1-u}=\sum_{i=0}^\infty u^i$.
We will show the same for $\mathcal{T}^{em}_\HH$.

A key Lemma for \cite{BKP} (and also for us) is the following.
\begin{lemma}\cite[Lemma 10.3]{BKP}\label{lem:dist_is_corank_plus_nullity}
For any $c\in \mathbb{Z}^E$ there is at least one $h\in H$ such that 
\begin{equation*}\label{eq:simultaneous_minimizer}
	d_1^<(H,c)=d_1^<(h,c)\text{ and }d_1^>(H,c)=d_1^>(h,c).
\end{equation*}
Hence $d_1(H,c)=d_1(h,c)=d_1^<(h,c)+d_1^>(h,c)$, and for any $h'$ with $d_1(H,c)=d_1(h',c)$, we have $d_1^<(H,c)=d_1^<(h',c)$ and $d_1^>(H,c)=d_1^>(h',c)$.
\end{lemma}

\begin{proof}[Proof of Theorem \ref{thm:embedding_Tutte_well_defined}]
	We show that $\tilde{\mathcal{T}}_\HH(u,v) = \mathcal{T}^{em}_\HH\left(\frac{1}{1-u}, \frac{1}{1-v} \right)$ as formal power series, where $\frac{1}{1-u}=\sum_{i=0}^\infty u^i$. As by \cite{BKP}, $ \tilde{\mathcal{T}}_\HH(u,v) = \mathcal{T}_\HH\left(\frac{1}{1-u}, \frac{1}{1-v} \right)$, this means that $\mathcal{T}^{em}_\HH$ does not depend on the chosen ribbon structure and basis, and it agrees with $\mathcal{T}_\HH$.

    By the definition of $\tilde{\mathcal{T}}_\HH$ and Theorem \ref{thm:Crapo_intervals_partition} \eqref{st:partition}, we have 

    \begin{align*}
	\tilde{\mathcal{T}}_\HH(u,v)=\sum_{c\in\mathbb{Z}^E} u^{d_1^>(H,c)}v^{d_1^<(H,c)}=
	\sum_{h\in H(\HH)}\sum_{c\in C_{em}(h)} u^{d_1^>(H,c)}v^{d_1^<(H,c)} \end{align*}
    By Theorem \ref{thm:Crapo_intervals_partition} \eqref{st:min} and Lemma \ref{lem:dist_is_corank_plus_nullity} this further equals
    \begin{align*}
    	=\sum_{h\in H(\HH)}\sum_{c\in C_{em}(h)} u^{d_1^>(h,c)}v^{d_1^<(h,c)}
   \end{align*}
   Now, by the definition of $C_{em}(h)$, we have
   \begin{align*} 
   	=\sum_{h\in H(\HH)} \left(\frac{1}{1-u}\right)^{oi_{em}(h)} \left(\frac{1}{1-v}\right)^{oe_{em}(h)} \left(\frac{1}{1-u}+\frac{1}{1-v}-1\right)^{ie_{em}(h)}=\\ 
    	=\mathcal{T}^{em}_\HH\left(\frac{1}{1-u}, \frac{1}{1-v} \right)
    \end{align*}
where the equalities are understood as equalities of formal power series, with $\frac{1}{1-u}=\sum_{i=0}^\infty u^i$.  
\end{proof}

\section{Proving the existence of the embedding Crapo decomposition}\label{sec:existence_of_Crapo}

This section is dedicated to proving Theorem \ref{thm:Crapo_intervals_partition}, the existence of a Crapo decomposition for embedding activities of hypergraphs. 

\subsection{Lemmas on Jaeger trees and embedding orders}

We collect some properties of Jaeger trees that we will need.
Throughout this section, fix a ribbon structure and basis for $\bip\HH$.

We first recall an alternative characterization of Jaeger trees from \cite{semibalanced} that uses an ordering of all spanning trees of $\bip \HH$. 

For any two trees $T$ and $T'$, their tours (Definition \ref{def:tour_of_a_tree})) have some initial segment that coincides (this segment might be empty), and in this initial segment the next node-edge pair is chosen the same way for the two trees. Hence, the first difference between the tours of $T$ and $T'$ needs to be that some node-edge pair $(x,xy)$ is treated differently for the two trees: $xy\in T-T'$ or $xy\in T'-T$.
This enables us to define an ordering $\prec$ of spanning trees of $\bip \HH$.
\begin{defn}[$\prec$]For two spanning trees $T$ and $T'$ of $\bip\HH$, we set $T\prec T'$ if for the first difference $(x, xy)$ of their respective tour, either $x$ is an emerald node and $xy\in T'-T$ or $x$ is a violet node and $xy\in T-T'$.
\end{defn}

The following is the key property of Jaeger trees.
\begin{thm}\cite[Theorems 7.1 and 7.6]{semibalanced}\label{thm:Jaeger_tree_is_the_first_realizer}
	For each hypertree $h$, there is a unique realizing hypertree $T$. Moreover, $T$ is the minimal tree according to $\prec$ among all trees representing $h$. 
\end{thm}

We note that the minimality in the above theorem implies that for a given hypertree, the unique realizing Jaeger tree can be computed greedily. This is explained in more detail in \cite[Section 7]{semibalanced}.

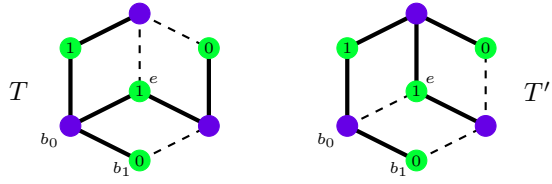
\begin{figure}[t]
	\begin{tikzpicture}[scale=.23]
		\begin{scope}[shift={(-8, 0)}]
			\draw [ultra thick] (18, 8.5) -- (14,6.5);
			\draw [ultra thick] (14, 6.5) -- (14,2);
			\draw [dashed, thick] (18, 8.5) -- (22,6.5);
			\draw [dashed, thick] (18, 8.5) -- (18,4);
			\draw [ultra thick] (14, 2) -- (18,4);
			\draw [ultra thick] (14, 2) -- (18,0);
			\draw [dashed, thick] (18, 0) -- (22,2);
			\draw [ultra thick] (22, 2) -- (18,4);
			\draw [ultra thick] (22, 2) -- (22,6.5);
			\draw [fill=e,e] (18, 0) circle [radius=0.6];
			\draw [fill=e,e] (18, 4) circle [radius=0.6];
			\draw [fill=e,e] (14, 6.5) circle [radius=0.6];
			\draw [fill=e,e] (22, 6.5) circle [radius=0.6];
			\draw [fill=v,v] (14, 2) circle [radius=0.6];
			\draw [fill=v,v] (22, 2) circle [radius=0.6];
			\draw [fill=v,v] (18, 8.5) circle [radius=0.6];
			
			\node at (14, 6.5) {{\tiny{$1$}}};
			\node at (22, 6.5) {{\tiny{$0$}}};
			\node at (18, 4) {{\tiny{$1$}}};
			\node at (18, 0) {{\tiny{$0$}}};
			
			\node at (11, 4) {{{$T$}}};
			\node at (12.8, 1.2) {{\tiny{$b_0$}}};
			\node at (17, -0.5) {{\tiny{$b_1$}}};
			\node at (18.8, 4.7) {{\tiny{$e$}}};
			
		\end{scope}
				
		\begin{scope}[shift={(8, 0)}]
				\draw [ultra thick] (18, 8.5) -- (14,6.5);
				\draw [ultra thick] (14, 6.5) -- (14,2);
				\draw [ultra thick] (18, 8.5) -- (22,6.5);
				\draw [ultra thick] (18, 8.5) -- (18,4);
				\draw [dashed,thick] (14, 2) -- (18,4);
				\draw [ultra thick] (14, 2) -- (18,0);
				\draw [dashed, thick] (18, 0) -- (22,2);
				\draw [ultra thick] (22, 2) -- (18,4);
				\draw [dashed, thick] (22, 2) -- (22,6.5);
				\draw [fill=e,e] (18, 0) circle [radius=0.6];
				\draw [fill=e,e] (18, 4) circle [radius=0.6];
				\draw [fill=e,e] (14, 6.5) circle [radius=0.6];
				\draw [fill=e,e] (22, 6.5) circle [radius=0.6];
				\draw [fill=v,v] (14, 2) circle [radius=0.6];
				\draw [fill=v,v] (22, 2) circle [radius=0.6];
				\draw [fill=v,v] (18, 8.5) circle [radius=0.6];
				
				\node at (14, 6.5) {{\tiny{$1$}}};
				\node at (22, 6.5) {{\tiny{$0$}}};
				\node at (18, 4) {{\tiny{$1$}}};
				\node at (18, 0) {{\tiny{$0$}}};
				
				\node at (25, 4) {{{$T'$}}};
				\node at (12.8, 1.2) {{\tiny{$b_0$}}};
				\node at (17, -0.5) {{\tiny{$b_1$}}};
				\node at (18.8, 4.7) {{\tiny{$e$}}};		
		\end{scope}
	\end{tikzpicture}
	\caption{Illustration for Example \ref{ex:tree_ordering}.}
	\label{fig:tree_ordering}
\end{figure}

\begin{ex}\label{ex:tree_ordering}
	Take Figure \ref{fig:tree_ordering}. Let the ribbon structure be the one induced by the positive orientation of the plane, and take the basis $(b_0,b_0b_1)$ according to the figure. 	
	$T$ (on the left panel) and $T'$ (on the right panel) are two spanning trees realizing the same hypertree. $T$ is Jaeger, while $T'$ is not, since $b_0e\notin T'$ is such that $(b_0,b_0e)$ precedes $(e,b_0e)$ in the tour of $T'$. 
	The first difference between the tours os $T$ and $T'$ is at $(b_0,b_0e)$, where $b_0$ is violet and $b_0b_1\in T-T'$. Hence indeed $T\prec T'$, cf. Theorem \ref{thm:Jaeger_tree_is_the_first_realizer}.
\end{ex}

For a spanning tree $T$ of $\bip\HH$ and edge $\varepsilon\in T$, we call the component of $T-\varepsilon$ containing $b_0$ the \emph{base component}.
The following property of fundamental cuts of Jaeger trees (an updated version of \cite[Lemma 6.13]{hyperBernardi}) makes them very useful to us. 

\begin{lemma}
	\label{l:char_Jaeger_cuts}\cite[Lemma 6.13]{hyperBernardi}
	Let $T$ be a Jaeger tree of $\bip\HH$, and $ev \in T$. 
	If $e'v'\in C^*(T,ev)-T$ is an edge such that its emerald endpoint $e'$ is in the base component, then $(e',e'v')$ comes earlier in the tour of $T$ than $(e,ev)$. 
\end{lemma}

\begin{lemma}\label{lem:emerald_nodes_of_E1_come_later}
	Let $T$ be a Jaeger tree of $\bip\HH$, realizing hypertree $h$, and let $ev \in T$ such that $e$ is in the base component of $T-ev$. Let $E_1$ be the set of emerald nodes not in the base component of $T-ev$. Then for any $f\in E_1$, we have $f>_h e$.
\end{lemma}
\begin{proof}
	The tour of $T$ only reaches the nodes of $E_1$ after traversing $ev$.
\end{proof}

\subsection{Lemmas on exchanges}

The following lemma is a special case of a well-known lemma in matroid theory, see for example \cite[Theorem 39.13]{Schrijver_CombOpt}. We sketch its proof for completeness. 
\begin{lemma}\cite[Theorem 39.13]{Schrijver_CombOpt} \label{lem:simultaneous_swap_in_trees}
	Let $\eta_1, \dots \eta_t, \varepsilon_1, \dots, \varepsilon_t$ be edges in the graph $G$, and let $T$ be a spanning tree such that $\eta_1, \dots \eta_t \in T$, $\varepsilon_1, \dots \varepsilon_t\notin T$, $\eta_i\in C(T,\varepsilon_i)$ $\forall i$ and $\eta_i\notin C(T,\varepsilon_j)$ for $j<i$. Then $T-\{\eta_1, \dots \eta_t\} + \{\varepsilon_1, \dots \varepsilon_t\}$ is also a spanning tree.
\end{lemma}
\begin{proof}
	Let $T_t=T- \eta_t + \varepsilon_t$. Then $T_t$ is a spanning tree since $\eta_t\in C(T,\varepsilon_t)$. Moreover, since $\eta_t\notin C(T,\varepsilon_j)$ for $j< t$, we have $C(T_t,\varepsilon_j)=C(T,\varepsilon_j)$ for all $j< t$. This implies that $T_{t-1}=T_t -\eta_{t-1} + \varepsilon_{t-1}$ is also a spanning tree, moreover, since $\eta_t, \eta_{t-1}\notin C(T,\varepsilon_j)$ for $j<t-1$, $C(T_{t-1},\varepsilon_j)=C(T,\varepsilon_j)$ for all $j< t-1$. We may continue switching edges like this until we arrive at $T'$.
\end{proof}

Let us state two variations of Lemma \ref{lem:simultaneous_swap_in_trees}, since these will be the formulations convenient to us.

\begin{lemma}\label{lem:partial_swap}
	Let $e_0v_0, \dots, e_{t-1}v_{t-1}$ and $e_1u_1,  \dots, e_tu_t$ be edges in the graph $G$, and let $T$ be a spanning tree such that $e_0v_0, \dots, e_{t-1}v_{t-1}\in T$, $e_1u_1, \dots, e_{t}u_{t} \notin T$ and $e_iv_i\in C(T,e_{i+1}u_{i+1})$ for $i=0, \dots t-1$. Then there exist indices $i_0=0 < i_1 < \dots < i_s = t$ such that
	$$T-\{e_{i_0}v_{i_0}, \dots, e_{i_{s-1}}v_{i_{s-1}}\} + \{e_{i_{1}}u_{i_{1}}, \dots, e_{i_s}u_{i_s}\}$$ is a spanning tree.
\end{lemma}
\begin{proof}
	Take the following auxiliary digraph $D$: Let the vertex set be $\{0, \dots, t\}$. Draw an edge from $i$ to $j$ if $e_iv_i \in C(T,e_ju_j)$. Then, by our assumptions, an edge points from $i$ to $i+1$ for each $0\leq i\leq t-1$, and there might also be other edges in the digraph. In any case, there is a directed path in $D$ from 0 to $t$. Choose a shortest directed path from $0$ to $t$ in $D$, and let the vertices of this path be $i_0=0, i_1, \dots , i_s=t$. Then $e_{i_{j}}v_{i_j} \in C(T,e_{i_{j+1}}u_{i_{j+1}})$ for each $0\leq j < s$ by construction, and $e_{i_{j}}v_{i_j} \notin C(T,e_{i_{j+1+r}}u_{i_{j+1+r}})$ for $r > 0$ because otherwise we could pick a shorter path from $0$ to $t$ in $D$. 
	
	Now we can apply Lemma \ref{lem:simultaneous_swap_in_trees} for $T$ with $\eta_j=e_{i_{s-j}}v_{i_{s-j}}$ and $\varepsilon_j=e_{i_{s-j+1}}u_{i_{s-j+1}}$. We get that
	$T-\{e_{i_0}v_{i_0}, \dots, e_{i_{s-1}}v_{i_{s-1}}\} + \{e_{i_{1}}u_{i_{1}}, \dots, e_{i_s}u_{i_s}\}$ is indeed a spanning tree. 
\end{proof}

\begin{lemma}\label{lemma:partial_swap_masik}
	Let $e_0v_0, \dots, e_{t-1}v_{t-1}$ and $e_1u_1,  \dots, e_tu_t$ be edges in the graph $G$, and let $T$ be a spanning tree such that $e_0v_0, \dots, e_{t-1}v_{t-1}\notin T$, $e_1u_1, \dots, e_{t}u_{t} \in T$ and $e_{i+1}u_{i+1}\in C(T,e_{i}v_{i})$ for $i=0, \dots t-1$. Then there exist 
	indices $i_0=0 < i_1 < \dots < i_s = t$ such that 	$$T+\{e_{i_0}v_{i_0}, \dots, e_{i_{s-1}}v_{i_{s-1}}\} - \{e_{i_{1}}u_{i_{1}}, \dots, e_{i_s}u_{i_s}\}$$ is a spanning tree.
\end{lemma}
\begin{proof}
	Let us again define an auxiliary digraph $D$ on the vertex set $\{0, \dots, t\}$. Draw an edge from $i$ to $j$ if $e_ju_j \in C(T,e_iv_i)$. Then, by our assumptions, an edge points from $i$ to $i+1$ for each $0\leq i\leq t-1$, and there might also be other edges in the digraph. In any case, there is a directed path in $D$ from 0 to $t$. Choose a shortest directed path from $0$ to $t$ in $D$, and let the vertices of this path be $i_0=0, i_1, \dots , i_s=t$. Then $e_{i_{j+1}}u_{i_{j+1}} \in C(T,e_{i_j}v_{i_j})$ for each $0\leq j < s$, and $e_{i_{j+1+r}u_{i_{j+1+r}}} \notin C(T,e_{i_j}v_{i_j})$ for $r > 0$ because otherwise we could pick a shorter path from $0$ to $t$ in $D$. 
 Now we can apply Lemma \ref{lem:simultaneous_swap_in_trees} for $T$ with $\eta_j=e_{i_{k+j}}u_{i_{k+j}}$ and $\varepsilon_j=e_{i_{k+j-1}}v_{i_{k+j-1}}$. We get that
	$T+\{e_{i_0}v_{i_0}, \dots, e_{i_{s-1}}v_{i_{s-1}}\} - \{e_{i_{1}}u_{i_{1}}, \dots, e_{i_s}u_{i_s}\}$ is indeed a spanning tree.
\end{proof}

\subsection{Three technical lemmas on first differences of Jaeger trees}

In this section, we collect the 3 main technical lemmas needed to prove Thorem \ref{thm:Crapo_intervals_partition}.

\begin{lemma}\label{lem:key}
	Let $h$ and $h'$ be hypertrees with respective Jaeger trees $T$ and $T'$. Suppose that the first difference between the tours of $T$ and $T'$ is at $(e,\varepsilon)$, where $\varepsilon\in T'-T$.
	
	Then there exist a hyperedge $f$ such that $h(f)>h'(f)$, moreover, $h-\mathbf{1}_f+\mathbf{1}_e$ and $h'+\mathbf{1}_f -\mathbf{1}_e$ are both hypertrees, and $e <_h f$ and $e <_{h'} f$.
\end{lemma}
\begin{proof}
	Recall that $<_h$ is defined as the order in which the emerald tour of $T$ reaches the nodes in $E$. 
	
	Take the fundamental cut $C^*(T', \varepsilon)$. We call the shore containing $b_0$ the \emph{base component} of $T'-\varepsilon$, moreover, we denote its set of violet nodes by $V_0$ and its set of emerald nodes by $E_0$. Also, we denote $V_1=V-V_0$ and $E_1=E-E_0$.
	
	Apply Lemma \ref{l:char_Jaeger_cuts} to the Jaeger tree $T'$. This implies that each edge of $C^*(T',\varepsilon)-e$ that has its emerald endpoint in the base component became current in the emerald tour of $T'$ before $(e, \varepsilon)$. Since the emerald tours of $T$ and $T'$ coincide until $(e, \varepsilon)$, these edges also became current in the tour of $T$, and as they are not in $T'$, they are also not in $T$.
	As also $\varepsilon\notin T$, we conclude:
	\begin{equation}\label{eq:T_beli_alapvagas}
		\text{$T$ has no edge from $C^*(T',\varepsilon)$ that has an endpoint in $E_0$.}
	\end{equation}
	
	We claim that for any $g\in E_1$, we have $e <_h g$ and $e<_{h'}g$. Indeed, by Lemma \ref{lem:emerald_nodes_of_E1_come_later},  we have $e <_{h'} g$. As the tours of $T$ and $T'$ agree until reaching $e$, this implies that $e$ is also reached before $g$ in $T$, hence also $e <_{h} g$. We will look for our hyperedge $f$ within $E_1$.
	
	Take the fundamental cycle $C(T,\varepsilon)$. A cycle and a cut always meet in an even number of edges, and $\varepsilon\in C^*(T',\varepsilon)\cap C(T,\varepsilon)$, hence the cycle $C(T,\varepsilon)$ contains at least one edge $\eta=e_1u_1$ of $T\cap C^*(T',\varepsilon)$.  
	For this $\eta$, we also have $\varepsilon\in C(T',\eta)$, as $\eta$ has one endpoint in $E_0\cup V_0$ and one endpoint in $E_1\cup V_1$, and the only edge of $T'$ connecting these two sets is $\varepsilon$. By \eqref{eq:T_beli_alapvagas}, $e_1$ cannot be in $E_0$, thus, $e_1\in E_1$.
	
	Since $\eta\in C(T,\varepsilon)$, the graph $T-\eta + \varepsilon$ is a spanning tree. Moreover, it realizes $h-\mathbf{1}_{e_1} + \mathbf{1}_e$ where $e <_h e_1$. 
	Also, $\varepsilon\in C(T',\eta)$, hence $T'-\varepsilon + \eta$ is a spanning tree, and it realizes $h' + \mathbf{1}_{e_1}-\mathbf{1}_e$ where $e <_{h'} e_1$.
	
	Thus, if $h(e_1)> h'(e_1)$, then $e_1$ is a suitable choice for $f$. 
	If $h(e_1)\leq h'(e_1)$, then we have to continue looking for an appropriate hyperedge $f$. We describe the general step of this process. 
	In a general step, we will have a spanning tree $\tilde{T}$ satisfying the following properties:
	\begin{equation}\label{eq:tildeT_conditions}
		\begin{split}
			\tilde{T}\subseteq T\cup T' \text{ is a spanning tree representing $h$,}\\
			\varepsilon\notin \tilde{T}.
		\end{split}
	\end{equation}
	At the beginning, we have $\tilde{T}=T$, that satisfies all requirements. In a general step, $\tilde{T}$ is typically not a Jaeger tree (as $T$ is the only Jaeger tree representing $h$). 
	However, since $\tilde{T}\subseteq T\cup T'$ and $\varepsilon\notin \tilde{T}$, we still have the following property:
	\begin{equation}\label{eq:tildeT_beli_alapvagas}
		\text{$\tilde{T}$ has no edge from $C^*(T',\varepsilon)$ that has an endpoint in $E_0$.}
	\end{equation} 
	
	Let $\varepsilon_0=\varepsilon$, and suppose that for some $t\geq 1$, we have a sequence of edges $\varepsilon_0, \eta_1, \varepsilon_1, \dots , \eta_{t-1}, \varepsilon_{t-1}, \eta_t$ with the following properties: 
	\begin{equation}\label{eq:edge_sequence_properties}
		\begin{split}
			\forall i\in[0, t-1]:\quad \varepsilon_i=e_iv_i\in T'-\tilde{T},\\ 
			\forall i\in[1, t]:\quad \eta_i=e_iu_i\in \tilde{T}-T',\\ 
			\forall i\in[1, t]:\quad \eta_{i+1}\in C^*(T',\varepsilon_i)\cap C(\tilde{T},\varepsilon_i)-\varepsilon_i,\\ 
			\forall i\in[0, t-1]:\quad \varepsilon_i\in C(T', \eta_{i+1}),\\ 
			\forall i\in [1,t]:\quad e_i\in E_1,\\ 
			\text{the edges $\varepsilon_0, \eta_1, \varepsilon_1, \dots \varepsilon_{t-1}, \eta_t$ are all distinct,}\\ 
			\text{$h'+\mathbf{1}_{e_t}-\mathbf{1}_e$ is a hypertree},\\
			\text{$h-\mathbf{1}_{e_t}+\mathbf{1}_e$ is a hypertree.}    
		\end{split}
	\end{equation}
	
	For $\tilde{T}=T$, we have just proved that such a set of edges can be chosen for $t=1$. Also, note that in that argument, we did not use the Jaeger property of $T$, only \eqref{eq:T_beli_alapvagas}. The analogue of that, \eqref{eq:tildeT_beli_alapvagas} is true for an arbitrary $\tilde{T}$. Hence for any $\tilde{T}$ satisfying \eqref{eq:tildeT_conditions}, we can choose $\eta_1$ such that \eqref{eq:edge_sequence_properties} holds (with $t=1$).
	
	Suppose that we have $\varepsilon_0, \eta_1, \varepsilon_1, \dots , \eta_{t-1}, \varepsilon_{t-1}, \eta_t$ for some $t\geq 1$, satisfying \eqref{eq:edge_sequence_properties}.
	If $h(e_t)> h'(e_t)$, then $f = e_t$ satisfies all properties required by the lemma, and the proof is complete.
	
	If $h(e_t)\leq h'(e_t)$, then we show that either we can add a new pair of edges $\varepsilon_{t}, \eta_{t+1}$ and still satisfy \eqref{eq:edge_sequence_properties}, or we can find a new tree $\tilde{T}$ satisfying \eqref{eq:tildeT_conditions} that has a smaller symmetric difference to $T'$ than the current $\tilde{T}$ (and start anew with this new $\tilde{T}$). If we prove this, then as $|T\Delta T'|$ is finite, moreover, $\{\varepsilon_0,\dots, \varepsilon_{t-1}\}$ is a subset of $T'-\tilde{T}\subseteq T'-T$ which is a finite set, the process cannot go on indefinitely. Hence there must be a moment when we conclude that $h(e_t)>h'(e_t)$, and this gives us an $f$ satisfying the requirements of the lemma, thereby finishing the proof.
	
	Hence let us suppose that $h(e_t)\leq h'(e_t)$.
	Note that $e_t\neq e$, since $e\in E_0$ but $e_t\in E_1$ by our assumption.
	There might be more than one indices $i$ such that $e_i=e_t$, but in any case, since $e_t\neq e$, the number of edges $\eta_i$ incident to $e_t$ (from $\tilde{T}-T'$) is one larger than the number of edges $\varepsilon_i$ (from $T'-\tilde{T}$). By \eqref{eq:edge_sequence_properties}, these edges are all different. 
	Since $h(e_t) = d_{\tilde{T}}(e_t)\leq h'(e_t)=d_{T'}(e_t)$, the degree of $e_t$ in $\tilde{T}-T'$ is at most the degree of $e_t$ in $T'-\tilde{T}$. 
	Thus, there is at least one so far unchosen edge $\varepsilon_t=e_tv_t$ incident to $e_t$ that is in $T'-\tilde{T}$.
	
	Take $C^*(T',\varepsilon_t)\cap C(\tilde{T},\varepsilon_t)$. As $\varepsilon_t$ is in this intersection, there need to be at least one more edge $\eta_{t+1}\in C^*(T',\varepsilon_t)\cap C(\tilde{T},\varepsilon_t)$ that is necessarily in $\tilde{T}-T'$. Let $\eta_{t+1}=e_{t+1}u_{t+1}$. 
	
	This $\eta_{t+1}$ might or might not agree with some previously chosen $\eta_i$, but before addressing this issue, let us note two things.
	
	Note that $\varepsilon_t\in C(T', \eta_{t+1})$ as $\varepsilon_t$ is the only edge of $T'$ connecting the two shores of $C^*(T',\varepsilon_t)$.
	
	Also, let us show that $e_{t+1}\in E_1$.
	Recall that $e_t\in E_1$. As $(E_0\cup V_0, E_1\cup V_1)$ is a fundamental cut of $T'$, one of the components of $T' - \varepsilon_t$ only contains vertices from $E_1\cup V_1$. 
	Hence $e_{t+1}\in E_0$ would imply that $\eta_{t+1}$ has endpoints in $E_0$ and $V_1$, thus, $\eta_{t+1}\in \tilde{T}\cap C^*(T',\varepsilon)$. However, by \eqref{eq:tildeT_beli_alapvagas} such an edge cannot have an endpoint in $E_0$, hence we conclude that $e_{t+1}\in E_1$. 
	
	From here, we distinguish 2 cases: 
	
	Case 1: $\eta_{t+1}$ does not agree with any of the edges $\eta_1, \dots , \eta_t$. In this case we show that $\varepsilon_0, \eta_1, \varepsilon_1, \dots , \eta_{t-1}, \varepsilon_{t}, \eta_{t+1}$ satisfies \eqref{eq:edge_sequence_properties}. To prove this, it is left to show that $h'+\mathbf{1}_{e_{t+1}}-\mathbf{1}_e$ is a hypertree, and $h-\mathbf{1}_{e_{t+1}}+\mathbf{1}_e$ is a hypertree. 
	
	We know that $\varepsilon_i\in C(T', \eta_{i+1})$, for each $0\leq i \leq t$. Apply Lemma \ref{lem:partial_swap} for $T'$, and $\varepsilon_0, \eta_1, \dots \varepsilon_t, \eta_{t+1}$. By the Lemma, $h'+\mathbf{1}_{e_{t+1}}-\mathbf{1}_e$ is indeed a hypertree. 
	 
	For each $0\leq i \leq t$ we also have $\eta_{i+1} \in C(\tilde{T},\varepsilon_{i})$. Hence we can apply Lemma \ref{lemma:partial_swap_masik} for $\tilde{T}$ and the edges $\varepsilon_0, \eta_1, \dots, \varepsilon_t, \eta_{t+1}$, and obtain that $h-\mathbf{1}_{e_{t+1}}+\mathbf{1}_e$ is a hypertree. 
	
	Case 2: $\eta_{t+1}=\eta_i$ for some $1\leq i \leq t$. In this case we show that we can choose a new $\tilde{T}$ satisfying \eqref{eq:tildeT_conditions} that has a strictly smaller symmetric difference to $T'$ than the current one.
	
	Note that for each $0\leq i \leq t$, $\eta_{i+1} \in C(\tilde{T},\varepsilon_{i})$. Hence we can apply Lemma \ref{lemma:partial_swap_masik} for $\tilde{T}$ and the edges $\varepsilon_i, \eta_{i+1}, \dots, \varepsilon_t, \eta_{t+1}=\eta_i$. The Lemma gives us a new tree $\tilde{T}_1=\tilde{T}+\varepsilon_{i_0}-\eta_{i_1}+ \varepsilon_{i_1} - \eta_{i_2} +\dots + \varepsilon_{i_{s-1}} - \eta_{i_s}$ where $\varepsilon_{i_0}=\varepsilon_{i}$ and $\eta_{i_s}=\eta_{t+1}=\eta_i$. Hence $\tilde{T}_1$ is a spanning tree that again realizes $h$, moreover, $\tilde{T}_1\subseteq \tilde{T}\cup T'\subseteq T\cup T'$. As neither of the $\eta_{i+1}, \dots, \varepsilon_t, \eta_{t+1}$ are incident to $e$, we also have $\varepsilon\notin \tilde{T}_1$. Also, since $\tilde{T}_1 \Delta T' \subseteq \tilde{T}\Delta T' - \varepsilon_i$, by choosing $\tilde{T}_1$ as our new $\tilde{T}$, we obtain a new tree satisfying \eqref{eq:tildeT_conditions}, that has a strictly smaller symmetric difference to $T'$ then the previous one. 
\end{proof}

\begin{lemma}\label{lemma:there_is_transfer_using the first difference}
	Let $J$ and $J'$ be Jaeger trees realizing hypertrees $x$ and $x'=x+\mathbf{1}_a -\mathbf{1}_b$. Suppose that the first difference between the tours of $J$ and $J'$ is at $(g,gv)$ with $gv\in J'-J$. Then 
	\begin{enumerate}[label=(\roman*)]
		\item\label{g_nem_b} $g\neq b$,
		\item\label{x+a-g} $x+\mathbf{1}_a -\mathbf{1}_g$ is a hypertree,
		\item\label{x+g-b} $x+\mathbf{1}_g -\mathbf{1}_b$ is a hypertree.
	\end{enumerate}
\end{lemma}
Note that $a=g$ is possible, and in that case $(ii)$ and $(iii)$ are true by definition.

\begin{proof}
	Apply Lemma \ref{lem:key} with $h=x$, $h'=x'$, $T=J$, $T'=J'$, and $e=g$. We get that there exist some hyperedge $f$ such that $x(f)>x'(f)$, moreover, $x-\mathbf{1}_f+\mathbf{1}_g$ and $x'+\mathbf{1}_f-\mathbf{1}_g$ are both hypertrees, and $g <_h f$, $g <_{h'} f$.
	
	As $x'=x+\mathbf{1}_a-\mathbf{1}_b$, the only hyperedge with $x(f)>x'(f)$ is $f=b$. Hence $x-\mathbf{1}_b+\mathbf{1}_g$ and $x'+\mathbf{1}_b-\mathbf{1}_g=x+\mathbf{1}_a-\mathbf{1}_g$ are both hypertrees. Finally, as $g <_h b$, we cannot have $g=b$. 
\end{proof}

\begin{lemma}\label{lemma:first_diff_for_Jaeger_of_a_transfer_is_earlier}
	Let $x$ and $x'=x+\mathbf{1}_a-\mathbf{1}_b$ be two hypertrees with respective Jaeger trees $J$ and $J'$. Suppose that the first difference in the tours of $J$ and $J'$ is at $(g, gv)$.  
	Then $g \leq_x a$, and $g\leq_x b$. 
\end{lemma}
Note that by the Jaeger property $g$ is an emerald node, hence $g \leq_x a$ and $g\leq_x b$ indeed makes sense.

\begin{prop}\label{prop:fundamental_cycle_of_a_Jaeger-tree}
	Suppose that $T$ is a Jaeger tree with $ev \notin T$, and suppose that the first node that we reach from $C(T,ev)$ in the tour of $T$ is not $e$. 
	Then, for each emerald node $f$ of $C(T,\varepsilon)$, 
	either $f=e$ or the last time $f$ is visited in the tour of $T$ is later than the last time $e$ is visited in the tour of $T$.
\end{prop}
\begin{proof}
	Let $g$ be the first node (either violet or emerald) of $C(T,ev)$ that is reached in the tour of $T$. Then $C(T,ev)$ consists of a path between $g$ and $e$, a path between $x$ and $v$, and the edge $ev$. By the Jaeger property, $ev$ is first reached at $e$, hence the path between $g$ and $e$ is visited first. The nodes on this path are ancestors of $e$ (or agree with $e$), hence their last visit is after the last visit of $e$. The rest of the vertices of $C(T,ev)$ (that is, the path between $g$ and $v$ excuding $x$) is only reached after the last visit of $e$, hence this is true also for their last visits.
\end{proof} 

\begin{proof}[Proof of Lemma \ref{lemma:first_diff_for_Jaeger_of_a_transfer_is_earlier}]
	It is enough to prove the statement if $gv\in J'-J$. Indeed, if we know the statement for $gv\in J'-J$, but we have $gv\in J-J'$, then we can apply the Lemma for $x'$ and $x=x'+\mathbf{1}_b-\mathbf{1}_a$, and obtain that $g\leq_{x'} b$ and $g\leq_{x'} a$. However, as the tours of $J$ and $J'$ agree until reaching $g$, $g\leq_{x'} b$ implies $g\leq_{x} b$ and $g\leq_{x'} a$ implies $g\leq_{x} a$.
	
	Hence let us suppose from now on that $gv\in J'-J$. 
	Apply Lemma \ref{lem:key} with $T=J$, $T'=J'$, $h=x$, $h'=x'$ and $e=g$. We get that there exist a hyperedge $f$ such that $g <_x f$ and $x(f)>x'(f)$. As $x'=x+\mathbf{1}_a-\mathbf{1}_b$, we can only have $f=b$, and hence $g<_x b$.
	
	Now let us prove $g \leq_x a$.
	Think of $J$ as a rooted tree, with root $b_0$ (the base vertex). Now we can talk about anchestors and descendants within $J$. Note that if an emerald vertex $e$ is a descendant of an emerald vertex $f$ in $T$, then we have $f <_x e$. 
	
	We start with showing that $g$ is not a descendant of $a$ in $J$. Suppose for a contradiction that $g$ is a descendant of $a$. This means that the last visit of $g$ in the tour of $J$ precedes the last visit of $a$. We will reach a contradiction by finding an alternative tree $J''$ realizing $x'$ such that $J''\prec J'$. This will be a contradiction by Theorem \ref{thm:Jaeger_tree_is_the_first_realizer} since we supposed that $J'$ is a Jaeger tree.
	
	Take $e_0=a$ and let $\varepsilon_0$ be an edge of $J'-J$ incident to $a=e_0$ (such an edge exist since $x'(a)>x(a)$). 
	Take $C(J,\varepsilon_0)$ and $C^*(J',\varepsilon_0)$. As the intersection of a cycle and a cut has an even number of edges, and $\varepsilon_0\in C(J,\varepsilon_0)\cap C^*(J',\varepsilon_0)$, there is also an edge 
	$\eta_1=e_1u_1\in C(J,\varepsilon_0)\cap C^*(J',\varepsilon_0)-J'$.
	
	We claim that the last visit of $e_1$ is after the last visit of $g$. We need to look at two cases. If $a$ is the first reached vertex within $C(J,\varepsilon_0)$, then $\varepsilon_0$ is the first considered edge of $C(J,\varepsilon_0)$ in the tour of $J$. Since $gv$ is supposed to be a first difference between $J$ and $J'$, but also $\varepsilon_0\in J'-J$, we need to arrive to $g$ before $(a,\varepsilon_0)$. As we supposed that $g$ is a descendant of $a$, in this case the last visit of $g$ is also before $(a,\varepsilon_0)$. All the vertices on $C(J,\varepsilon_0)$ except for $a$ are reached after $(a,\varepsilon_0)$ hence we conclude that the last visit of $e_1$ is after the last visit of $g$.	 
	If $C(J,\varepsilon_0)$ is not reached at $a$, then Proposition \ref{prop:fundamental_cycle_of_a_Jaeger-tree} implies that the last visit of $e_1$ is after the last visit of $a$, and hence also after the last visit of $g$.
	
	If $e_1=b$, then $J''=J+\varepsilon_0 -\eta_1$ is a spanning tree realizing $h'$. Moreover, as $\varepsilon_{0}, \eta_1 \in J\Delta J'$, and the last visit of both $e_0$ and $e_1$ are after the last visit of $g$, we still have $gv\notin J''$, and the tour of $J''$ agrees with the tour of $J'$ until reaching $(g,gv)$. Hence $J''\prec J'$, contradicting Theorem \ref{thm:Jaeger_tree_is_the_first_realizer}.
	
	If $e_1\neq b$, we continue looking for a tree realizing $x'$. Let us describe the general step of this process. In a general step, we will have a tree $\tilde{J}\subset J\cup J'$ with $gv\in\tilde{J}$, realizing $x'$. We also have some $t\geq 0$ and edges $\varepsilon_0, \eta_1, \varepsilon_1,  \dots, \eta_{t-1}, \varepsilon_{t-1}, \eta_t$ (in case $t=0$, we do not have any edges) with the following properties: For $i\in [0,t-1], \varepsilon_i=e_iv_i\in \tilde{J}-J$, for $i\in [1,t], \eta_i=e_iu_i\in J-\tilde{J}$.
	The edges $\varepsilon_0, \eta_1, \varepsilon_1, \dots \varepsilon_{t-1}, \eta_t$ are all distinct, $\eta_{i+1}\in C(J,\varepsilon_i)$ for $i\in [0, t-1]$, and $\varepsilon_i\in C(\tilde{J}, \eta_{i+1})$ for $i\in [0,t-1]$. We moreover require that $e_0=a$, 
	and for each $1\leq i\leq t$, the last visit of $e_i$ in the tour of $T$ is after the last visit of $g$. 
	
	At the beginning, $\tilde{J}=J'$ and $t=0$ obviously satisfies the requirements.
	
	We show that if $e_{t} \neq b$ (which also includes the case if $t=0$), then we can either find an additional pair of edges $\varepsilon_t$ and $\eta_{t+1}$ such that the extended sequence of edges still satisfies the required properties, or we can find a new $\tilde{J}$ realizing $x'$ such that $|\tilde{J}\Delta J|$ strictly decreases (and take $t=0$).
	
	If $e_t\neq b$, then in case if $e_t\neq a$, we have $x(e_t)=x'(e_t)$ and the number of so far chosen edges from $J-\tilde{J}$ incident to $e_t$ is one larger than the number of edges from $\tilde{J}-J$. Hence in this case there is an edge $\varepsilon_t=e_tv_t\in \tilde{J}-J$ that we have not chosen so far.  
	If $e_t = a$, then $x'(e_t)=x(e_t)+1$, but the number of so far chosen edges from $J-\tilde{J}$ incident to $e_t$ is equal to the number of edges from $\tilde{J}-J$. Hence in this case, too, there is an edge $\varepsilon_t=e_tv_t\in \tilde{J}-J$ that we have not chosen so far. As the intersection of a cycle and a cut has an even number of edges, in both cases there is also an edge 
	$\eta_{t+1}=e_{t+1}u_{t+1}\in C(J,\varepsilon_t)\cap C^*(\tilde{J},\varepsilon_t)-\tilde{J}$.
	This also implies $\varepsilon_t\in C(\tilde{J}, \eta_{t+1})$.
	
	We claim that the last visiting time of $e_{t+1}$ is after the last visiting time of $g$. If $C(J,e_tv_t)$ was not reached at $e_t$ in the tour of $J$, then by Proposition \ref{prop:fundamental_cycle_of_a_Jaeger-tree}, the last visiting time of $e_{t+1}$ is after the last visiting time of $e_t$, hence it is after the last visiting time of $g$.
	If $C(J,e_tv_t)$ was reached at $e_t$ in the tour of $J$, then $e_{t+1}$ is a descendant of $e_t$. There are two cases: either $g$ is a descendant of $e_t$ or the first reaching time of $e_t$ is after the last visiting time of $g$. In the latter case, it is obvious that the last visiting time of $e_{t+1}$ is after the last visiting time of $g$. If $e_t$ is an ancestor of $g$, then as $e_tv_t\in \tilde{J}-J$, the pair $(g,gv)$ needs to precede $(e_t, e_tv_t)$ in the tour of $J$, and hence by the time the tour of $J$ is in $(e_t, e_tv_t)$, the node $g$ is finished. Also, $e_{t+1}$ is only reached after $(e_t, e_tv_t)$. Hence indeed, its last (and also first) visiting time is after the last visiting time of $g$.
	
	If $\eta_{t+1}$ does not agree with $\eta_i$ for any $i\in [1,t]$, then we extended our sequence of edges such that $\varepsilon_0, \eta_1, \varepsilon_1,  \dots, \varepsilon_{t}, \eta_{t+1}$ satisfies all required properties.
	
	If $\eta_{t+1}=\eta_i$ for some previous $i$, then  apply Lemma \ref{lemma:partial_swap_masik} to $\tilde{J}$ and the sequence of edges $\varepsilon_{t+1}, \eta_{t}, \varepsilon_{t}, \dots, \varepsilon_{i+1}, \eta_i$. We get that there is a sequence of edges $t+1=i_0 > i_1 > \dots > i_s=i$ such that $\tilde{J}_1=\tilde{J}+\{e_{i_0}v_{i_0}, \dots, e_{i_{s-1}}v_{i_{s-1}}\} - \{e_{i_{1}}u_{i_{1}}, \dots, e_{i_s}u_{i_s}\}$ is a spanning tree. As $e_{i_0}=e_{t+1}=e_i=e_{i_s}$, $\tilde{J}_1$ again realizes $x'$. Moreover, we have $\tilde{J}_1 \subseteq \tilde{J}\cup J \subseteq J' \cup J$, and as the last visit of each $e_j$ ($0\leq j \leq t+1$) is after te last visit of $g$, the edge $gv$ was not among the edges that we switched, hence also $gv\in \tilde{J}_1$. As $e_iu_i \in \tilde{J}-J$, we have $|\tilde{J}_1\Delta J|<|\tilde{J}\Delta J|$, hence we can choose $\tilde{J}_1$ as our new $\tilde{J}$.
	
	Altogether, as $J'\Delta J$ is finite, we cannot have the above two cases infinitely many times, hence after a while, we need to have a case where $e_t=b$.
	
	By Lemma \ref{lemma:partial_swap_masik} applied to $J$ and our sequence of edges, there exist a subset of these edges such that $i_0=0 < i_1 < \dots < i_s=t$ and 
	$J''=J + e_{i_0}v_{i_0} - e_{i_1}u_{i_1} + e_{i_1}v_{i_1} - \dots + e_{i_{s-1}}v_{i_{s-1}}- e_{i_s}u_{i_s}$ is a spanning tree.
	$J''$ realizes $x'=x+\mathbf{1}_a-\mathbf{1}_b$. Moreover, as 
	for each $i$, the last visit of $e_i$ is after the last visit of $g$, none of the $e_i$ agree with $g$. Hence none of the edges $e_0v_0, e_{1}u_{1}, e_{1}v_{1}, \dots, e_{t-1}v_{t-1}, e_{t}u_{t}$ agrees with $gv$. As $(g,gv)$ is the first difference between $J$ and $J'$, and $J''\subseteq J\cup \tilde{J}\subseteq J \Delta J'$, the first difference between $J$ and $J''$ is at $(g,gv)$, where $gv\in J'-J''$. This means that $J''\prec J'$, contradicting Theorem \ref{thm:Jaeger_tree_is_the_first_realizer}, as both trees realize $x'$. This proves that $a$ cannot be an ancestor of $g$ in $J$.
	
	If $a$ is not an ancestor of $g$, then there are three possibilities. If $g$ is an ancestor of $a$ or $g=a$, then $g \leq_x a$.
	The only other possibility is that neither of $a$ and $g$ are ancestors of each other. In this case, the last visit of the earlier visited node is before the first visit of the other node. As the first difference between the tours of $J$ and $J'$ is at $g$, and there is also at least one edge $\varepsilon_0$ of $J'-J$ incident to $a$, we conclude that $g$ needs to be visited earlier than $a$, otherwise $\varepsilon_0$ would be visited before $gv$. Hence we also conclude that $g\leq_x a$. 
\end{proof}

\subsection{Proof of Theorem \ref{thm:Crapo_intervals_partition}}

\begin{thm}\label{thm:separation}
    Fix an arbitrary ribbon structure and basis.
	If $h$ and $h'$ are two hypertrees of a hypergraph $\mathcal{H}$, then $C_{em}(h)\cap C_{em}(h')=\emptyset$. 
\end{thm}
\begin{proof}
Let $T$ and $T'$ be the respective Jaeger trees of $h$ and $h'$. By symmetry, we may suppose that for the first difference $(e,\varepsilon)$ between the tours of $T$ and $T'$, we have $\varepsilon\in T-T'$. (As both $T$ and $T'$ are Jaeger trees, $e$ needs to be a hyperedge.)

Apply Lemma \ref{lem:key}. Notice that the hyperedge $f$ given by the Lemma has the following properties:
 \begin{equation}\label{eq:tilde_e_properties}
 \begin{split}
 h(f)>h'(f), \\ 
 f\text{ is not internally embedding active in $h$,}\\ 
 f\text{ is not externally embedding active in $h'$}.    
 \end{split}
 \end{equation}
  
This implies that for each $x \in C_{em}(h)$ and $y\in C_{em}(h')$ we have $x(f) \geq h(f) > h'(f) \geq y(f)$, hence $x$ and $y$ cannot be the same.
\end{proof}

\begin{thm}\label{thm:each_point_in_some_interval}
	For any $c\in \mathbb{Z}^E$, there is at least one hypertree $h$ such that $c\in C_{em}(h)$, moreover, $d_1(h,c)=d_1(H,c)$.
\end{thm}
Before proving the theorem, let us recall some constructions, and prove some lemmas. 

For a given vector $c\in\mathbb{Z}^E$, let $H_c=\{h\in H \mid d_1(H,c)=d_1(h,c)\}$. We are looking for a hypertree $h\in H_c$ such that $c\in C_{em}(h)$.

Our method will be the following: We define a complete ordering $<_c$ of the hypertrees, that orders them according to the first differences between their Jaeger trees in a way depending on $c$. Then, we show that if a hypertree $h\in H_c$ does not have $c\in C_{em}(h)$, then we can find another hypertree $h'\in H_c$ with $h' <_c h$. As there are finitely many hypertrees, this implies that at some point, we need to find a hypertree $h\in H_c$ with $c\in C_{em}(h)$.

To define the ordering $<_c$, we need some more preparations.

\begin{prop}\label{prop:external_internal_partition_for_c}
If for an edge $e$, there exist $h\in H_c$ with $h(e)> c(e)$, then for each $h'\in H_c$, we have $h'(e)\geq c(e)$, and symmetrically, if there exist $h\in H_c$ with $h(e)< c(e)$, then for each $h'\in H_c$, we have $h'(e)\leq c(e)$.
\end{prop}

\begin{proof}
	Suppose for a contradiction that there is a hyperedge $e$ such that there are $h,h'\in H_c$ with $h(e)< c(e) < h'(e)$. Then, by Proposition \ref{prop:exchange_property_for_hypertrees},
	there exist a hyperedge $f$ with $h(f)>h'(f)$ such that $h+\mathbf{1}_{e}-\mathbf{1}_{f}$ and $h'-\mathbf{1}_{e}+\mathbf{1}_{f}$ are both hypertrees.
	
	Then $d_1(h+\mathbf{1}_{e}-\mathbf{1}_{f},c)\leq d_1(h,c)$. As $d_1(h,c)$ is minimal among all hypertrees, we need to have $d_1(h+\mathbf{1}_{e}-\mathbf{1}_{f},c)=d_1(h,c)$. This implies that $c(f)\geq h(f)$. Similarly, $d_1(h'-\mathbf{1}_{e}+\mathbf{1}_{f},c)\leq d_1(h',c)$. As $d_1(h',c)$ is also minimal among all hypertrees, we need to have $d_1(h'-\mathbf{1}_{e}+\mathbf{1}_{f},c)=d_1(h',c)$. Hence $c(f)\leq h'(f)$. These together imply that $h(f)\leq h'(f)$, which is a contradiction.
\end{proof}

We say that a hyperedge $e$ is \emph{external} with respect to $c$ if there is a hypertree $h\in H_c$ such that $h(e) < c(e)$. We say that a hyperedge $e$ is \emph{internal} with respect to $c$ if there exist a hypertree $h\in H_c$ with $h(e) > c(e)$. Then by Proposition \ref{prop:external_internal_partition_for_c}, each hyperedge is either internal or external, or neither of them, in which case each $h\in H_c$ has $h(e)=c(e)$. Moreover, for an internal hyperedge, each $h\in H_c$ has $h(e) \geq c(e)$, and for an external hyperedge, each $h\in H_c$ has $h(e) \leq c(e)$.

Now we are ready to define the ordering $<_c$.
\begin{defn}[$<_c$]
	Let $h$ and $h'$ be two hypertrees with respective Jaeger trees $T$ and $T'$. Suppose that the fist difference between the tours of $T$ and $T'$ is at $(e,ev)$, where $ev\in T-T'$. Then
	\begin{enumerate}
		\item{} if $e$ is an external hyperedge, we define $h' <_c h$
		\item{} if $e$ is an internal hyperedge, we define $h <_c h'$
		\item{} if $e$ is neither an external, nor an internal hyperedge, we define $h <_c h'$.
	\end{enumerate}
\end{defn}
\begin{prop}
	For any $c\in\mathbb{Z}^E$, the above defined $<_c$ is a complete ordering of $H$.
\end{prop}
\begin{proof}
	As any two trees are comparable, we only need to prove transitivity. Suppose that $T_1<_c T_2$ and $T_2<_c T_3$. Assume that the first difference between $T_1$ and $T_2$ is at $(e,\varepsilon)$ and the first difference between $T_2$ and $T_3$ is at $(e',\varepsilon')$.
	
	Suppose that $e$ and $e'$ are both external with respect to $c$. Then we might have $e=e'$, but in any case, $\varepsilon\neq \varepsilon'$ since $\varepsilon\in T_2$ while $\varepsilon'\notin T_2$.
	  
	If $(e,\varepsilon)$ precedes $(e',\varepsilon')$ in the tour of $T_2$, then the first difference between $T_1$ and $T_3$ is $(e,\varepsilon)$, and this orders them as $T_1$ and $T_2$ are ordered, i.e., $T_1\prec T_3$. If $(e',\varepsilon')$ precedes $(e,\varepsilon)$ in the tour of $T_2$, then the first difference between $T_1$ and $T_3$ is $(e',\varepsilon')$, and this implies that they are ordered as $T_2$ and $T_3$ are, that is, again, $T_1\prec T_3$.
	
	If $e$ is external and $e'$ is internal with respect to $c$, then $e\neq e'$ and thus $\varepsilon\neq \varepsilon'$ again. The rest of the argument can be repeated.
	
	The reasoning is analogous also in the cases when $e$ and $e'$ are both internal, and when $e$ is internal and $e'$ is external with respect to $c$.
\end{proof}

\begin{proof}[Proof of Theorem \ref{thm:each_point_in_some_interval}]
	Let us spell out the property $c\in C_{em}(h)$ in another way:
	\begin{equation}\label{eq:external}
	\begin{split}
	\text{If $e$ is an external edge of $c$ and $c(e)> h(e)$},\\ 
	\text{then $e$ is externally embedding active in $h$,}
	\end{split}
	\end{equation}
	\begin{equation}\label{eq:internal}
	\begin{split}
	\text{if $e$ is an internal edge of $c$ and $c(e)< h(e)$},\\ 
	\text{then $e$ is internally embedding active in $h$.}
	\end{split}
	\end{equation}
	Clearly, $c\in C(h)$ is equivalent to \eqref{eq:external} plus \eqref{eq:internal}.
		
	Take an arbitrary hypertree $h\in H_c$. If $h$ satisfies \eqref{eq:external} and \eqref{eq:internal} then we are done. If not, then we will modify $h$ to get a new hypertree $h'\in H_c$ with $h' <_c h$. As $<_c$ is a complete ordering on hypertrees and there are finitely many hypertrees, at some point we need to have a hypertree $h$ satisfying \eqref{eq:external} and \eqref{eq:internal}.
	
    If $h$ does not satisfy either \eqref{eq:external} or \eqref{eq:internal}, then there is either an external hyperedge $e$ not satisfying \eqref{eq:external} or an internal hyperedge $e$ not satisfying \eqref{eq:internal} (or both). We next show how to modify $h$ in these two cases.
	
	Suppose that $e$ is an external hyperedge such that $c(e)> h(e)$ but $e$ is not externally embedding active in $h$.
	By the definition of external embedding activity, this means that there exist a hyperedge $f <_h e$ such that $h'= h+\mathbf{1}_e - \mathbf{1}_f$ is also a hypertree. 
	Note that such an $f$ is necessarily an external hyperedge, and $h'\in H_c$. Indeed, $h\in H_c$ implies $d_1(h,c)\leq d_1(h',c)$. As $|c(e)-h'(e)|<|c(e)-h(e)|$ we conclude that $c(f)> h'(f)$, thus, $f$ is also external, and $h'\in H_c$. 
	
	We show that we can choose $f$ such that $h' <_c h$. More precisely, we show that if we choose $f$ to be the earliest emerald node according to $<_h$ such that $h+\mathbf{1}_e-\mathbf{1}_f$ is a hypertree, then for the obtained $h'=h+\mathbf{1}_e-\mathbf{1}_f$ and its Jaeger tree $T'$, the first difference between the tours of $T$ and $T'$ is at a pair $(f, fv)$, where $fv\in T-T'$. As $f$ is an external hyperedge (with respect to $c$), this means that $h' <_c h$. 
	
	Suppose for a contradiction that $f$ is chosen to be the earliest emerald node according to $<_h$ such that $h+\mathbf{1}_e-\mathbf{1}_f$ is a hypertree, but the above properties do not hold for the Jaeger trees $T$ and $T'$. There are 3 ways in which the property can be violated.
	
	Case 1: The first difference between $T$ and $T'$ is at $(g,gv)$ where $g\neq f$, and $gv\in T-T'$. 
	
	By Lemma \ref{lemma:first_diff_for_Jaeger_of_a_transfer_is_earlier} applied to $x=h$, $x'=h'$, $a=e$, $b=f$, in this case $g\leq_h f$ and $g\leq_h e$. 
	By part \ref{x+a-g} of Lemma \ref{lemma:there_is_transfer_using the first difference}, applied with $J=T'$, $J'=T$, $x=h'$, $a=f$, $b=e$,
	we get that $h'+\mathbf{1}_f -\mathbf{1}_g = h +\mathbf{1}_e -\mathbf{1}_g$ is also a hypertree. This contradicts the fact that $f$ was the earliest emerald node according to $<_h$ such that $h+\mathbf{1}_e-\mathbf{1}_f$ is a hypertree.
	
	Case 2: The first difference between $T$ and $T'$ is at $(g,gv)$ where $g\neq f$, and $gv\in T'-T$. 
	
	By Lemma \ref{lemma:first_diff_for_Jaeger_of_a_transfer_is_earlier} in this case $g\leq_h f$ and $g\leq_h e$. 
	By part \ref{x+a-g} of Lemma \ref{lemma:there_is_transfer_using the first difference}, applied for $J=T$, $J'=T'$, $x=h$, $a=e$, $b=f$,
	we get that $h+\mathbf{1}_e -\mathbf{1}_g$ is also a hypertree. This contradicts the fact that $f$ was the earliest emerald node according to $<_h$ such that $h+\mathbf{1}_e-\mathbf{1}_f$ is a hypertree.
	
	Case 3: The first difference between $T$ and $T'$ is at $(f,fv)$, but we have $fv\in T'-T$. This is impossible by \ref{g_nem_b} of Lemma \ref{lemma:there_is_transfer_using the first difference} (applied for $J=T$, $J'=T'$, $a=e$, $b=f$).
	
	Now let us look at the case if there is an internal hyperedge violating \eqref{eq:internal}.
	
	Let $e$ be an internal hyperedge such that $c(e) < h(e)$ and $e$ is not internally embedding active in $h$.
	By the definition of internal embedding activity, this means that there exist a hyperedge $f <_h e$ such that $h'= h - \mathbf{1}_e + \mathbf{1}_f$ is also a hypertree. 
	Note that $f$ is necessarily an internal hyperedge, and $h'\in H_c$. Indeed, $h\in H_c$ implies $d_1(h,c)\leq d_1(h',c)$. As $|c(e)-h'(e)|<|c(e)-h(e)|$ we conclude that $c(f)< h'(f)$, thus, $f$ is also internal, and $h'\in H_c$. 
	
    We show that we can choose $f$ such that $h' <_c h$. More precisely, we will show that if we choose $f$ to be the earliest emerald node according to $<_h$ such that $h-\mathbf{1}_e+\mathbf{1}_f$ is a hypertree, then for the obtained $h'=h-\mathbf{1}_e+\mathbf{1}_f$ and its Jaeger tree $T'$, the first difference between the tours of $T$ and $T'$ is at a pair $(f, fv)$, where $fv\in T'-T$. As $f$ is an internal hyperedge (with respect to $c$), this means that $h' <_c h$. 
    
    Suppose for a contradiction that $f$ is chosen to be the earliest emerald node according to $<_h$ such that $h-\mathbf{1}_e+\mathbf{1}_f$ is a hypertree, but the above properties do not hold for the Jaeger trees. There are 3 ways in which the property can be violated.
    
    Case 1: The first difference between $T$ and $T'$ is at $(g,gv)$ where $g\neq f$, and $gv\in T-T'$. 
    
    By Lemma \ref{lemma:first_diff_for_Jaeger_of_a_transfer_is_earlier} applied with $x=h$, $x'=h'$, $a=f$ and $b=e$, in this case $g\leq_h f$ and $g\leq_h e$. 
    
    By \ref{x+g-b} of Lemma \ref{lemma:there_is_transfer_using the first difference} applied with $J=T'$, $J'=T$, $x=h'$, $x'=h$, $a=e$ and $b=f$, we get that  
    $h'+\mathbf{1}_g -\mathbf{1}_f=h -\mathbf{1}_e + \mathbf{1}_g$ is a hypertree. This contradicts the fact that $f$ was the earliest emerald node according to $<_h$ such that $h-\mathbf{1}_e+\mathbf{1}_f$ is a hypertree.
    
    Case 2: The first difference between $T$ and $T'$ is at $(g,gv)$ where $g\neq f$, and $gv\in T'-T$. 
    
    By Lemma \ref{lemma:first_diff_for_Jaeger_of_a_transfer_is_earlier} in this case $g\leq_h f$ and $g\leq_h e$. 
    
    By \ref{x+g-b} of Lemma \ref{lemma:there_is_transfer_using the first difference} applied with $J=T$, $J'=T'$, $x=h$, $x'=h'$, $a=f$ and $b=e$, we get that $h+\mathbf{1}_g-\mathbf{1}_e$ is a hypertree. This again contradicts the fact that $f$ was the earliest emerald node according to $<_h$ such that $h-\mathbf{1}_e+\mathbf{1}_f$ is a hypertree.
    
    Case 3: The first difference between $T$ and $T'$ is at $(f,fv)$, but we have $fv\in T-T'$. This is impossible by \ref{g_nem_b} of Lemma \ref{lemma:there_is_transfer_using the first difference} applied for $J=T'$, $J'=T$, $x=h'$, $x'=h$, $a=e$ and $b=f$.
\end{proof}

\begin{proof}[Proof of Theorem \ref{thm:Crapo_intervals_partition}]
	The statement follows from Theorems \ref{thm:separation} and \ref{thm:each_point_in_some_interval}.
\end{proof}

\section{Open questions concerning variants of embedding activities}
\label{sec:open}

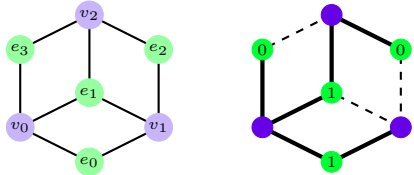
\begin{figure}[h]
	\begin{tikzpicture}[scale=.23]
			\begin{scope}[shift={(7, 0)}]
					\draw [dashed, thick] (18, 8.5) -- (14,6.5);
					\draw [ultra thick] (14, 6.5) -- (14,2);
					\draw [ultra thick] (18, 8.5) -- (22,6.5);
					\draw [ultra thick] (18, 8.5) -- (18,4);
					\draw [ultra thick] (14, 2) -- (18,4);
					\draw [ultra thick] (14, 2) -- (18,0);
					\draw [ultra thick] (18, 0) -- (22,2);
					\draw [dashed, thick] (22, 2) -- (18,4);
					\draw [dashed, thick] (22, 2) -- (22,6.5);
					\draw [fill=e,e] (18, 0) circle [radius=0.6];
					\draw [fill=e,e] (18, 4) circle [radius=0.6];
					\draw [fill=e,e] (14, 6.5) circle [radius=0.6];
					\draw [fill=e,e] (22, 6.5) circle [radius=0.6];
					\draw [fill=v,v] (14, 2) circle [radius=0.6];
					\draw [fill=v,v] (22, 2) circle [radius=0.6];
					\draw [fill=v,v] (18, 8.5) circle [radius=0.6];
					
					\node at (14, 6.5) {{\tiny{$0$}}};
					\node at (22, 6.5) {{\tiny{$0$}}};
					\node at (18, 4) {{\tiny{$1$}}};
					\node at (18, 0) {{\tiny{$1$}}};
					
				\end{scope}
			\begin{scope}[shift={(-7, 0)}]
					\draw [thick] (18, 8.5) -- (14,6.5);
					\draw [thick] (14, 6.5) -- (14,2);
					\draw [thick] (18, 8.5) -- (22,6.5);
					\draw [thick] (18, 8.5) -- (18,4);
					\draw [thick] (14, 2) -- (18,4);
					\draw [thick] (14, 2) -- (18,0);
					\draw [thick] (18, 0) -- (22,2);
					\draw [thick] (22, 2) -- (18,4);
					\draw [thick] (22, 2) -- (22,6.5);
					\draw [fill, le] (18, 0) circle [radius=0.8];
					\draw [fill, le] (18, 4) circle [radius=0.8];
					\draw [fill, le] (14, 6.5) circle [radius=0.8];
					\draw [fill, le] (22, 6.5) circle [radius=0.8];
					\draw [fill, lv] (14, 2) circle [radius=0.8];
					\draw [fill, lv] (22, 2) circle [radius=0.8];
					\draw [fill, lv] (18, 8.5) circle [radius=0.8];
					
					\node at (14, 6.5) {{\tiny{$e_3$}}};
					\node at (22, 6.5) {{\tiny{$e_2$}}};
					\node at (18, 4) {{\tiny{$e_1$}}};
					\node at (18, 0) {{\tiny{$e_0$}}};
					
					\node at (14, 2) {{\tiny{$v_0$}}};
					\node at (22, 2) {{\tiny{$v_1$}}};
					\node at (18, 8.5) {{\tiny{$v_2$}}};
				\end{scope}
		\end{tikzpicture}
	\caption{A violet Jaeger tree for the ribbon structure induced by the positive orientation of the plane, with basis $(v_0,v_0e_0)$.}\label{fig:violet_Jaeger}
\end{figure}

Let us recall an alternative definition for embedding activities from \cite{hyperBernardi}. In this paper, Jaeger trees are defined via the rule that each non-tree edge is first seen at its emerald endpoint. 
However, as mentioned in Remark \ref{rem:emerald_vs_violet_Jaeger_tree}, one could also define Jaeger trees with the rule that each non-tree edge is first seen at its violet endpoint (see Figure \ref{fig:violet_Jaeger} for an example). 

In \cite{hyperBernardi}, both types of trees are investigated, and are called emerald, and violet Jaeger trees, respectively. 

It is also true that for any hypertree $h$ of $\HH$, there is a unique violet Jaeger tree representing $h$. Hence one could potentially associate an ordering to $h$ using the violet Jaeger tree. There are in fact two natural ways to do this.

One possibility is to define the ordering $<^{violet}_h$ as the ordering in which emerald nodes are reached (that is, become current node) in the tour of $T$, where $T$ is the violet Jaeger tree representing $h$.

Another possibility is to take the ordering $<^{violet'}_h$ which is the ordering in which emerald nodes appear as an endpoint of the current edge in the tour of $T$, where $T$ is the violet Jaeger tree representing $h$.

\begin{ex}
	Figure \ref{fig:violet_Jaeger} shows the unique violet Jaeger tree $T$ for the hypertree $h$ with $h(e_0)=h(e_1)=1$, $h(e_2)=h(e_3)=0$. We have $e_0 <_h^{violet} e_1 <_h^{violet} e_2 <_h^{violet} e_3$ since this is the order in which emerald nodes become current in the tour of $T$. However, we have $e_0 <_h^{violet'} e_2 <_h^{violet'} e_1 <_h^{violet'} e_1 <_h^{violet'} e_3$, since $(v_1, v_1e_2)$ precedes $(v_1, v_1e_1)$ in the tour of $T$.
\end{ex}

Note that for an emerald Jaeger tree, these two orders coincide, since if an emerald endpoint $e$ first appears as an endpoint of the current edge such that the current node-edge pair is $(v,ve)$, then in the next step, the current node needs to be $e$. This is not true, however, if $T$ is a violet Jaeger tree, as also the previous example shows.

In some sense, $<^{violet}_h$ seems to be the more natural ordering, however, it is easy to come up with examples where it does not give $\mathcal{T}_\HH$.

However, $<^{violet'}_h$ seems to yield $\mathcal{T}_\HH$.

\begin{conj}
Embedding activities defined via $<^{violet'}_h$ produce the polynomial $\mathcal{T}_\HH$ for hypergraphs.
\end{conj}

\section{Courtiel's activities}\label{s:Delta_activities}

In \cite{Courtiel}, Courtiel gave a notion of activities for graphs, and more generally, matroids, that he calls $\Delta$-activities. $\Delta$-activities induce a Crapo-decomposition, and generalize several different notions of Tutte-descriptive activities, among others, activities with respect to a fixed ordering, and embedding activities (of graphs) \cite{Courtiel}.
He also gave a conjecture that any notion of activity for graphs that induces a Crapo-decomposition falls within the class of $\Delta$-activities. 

We show that Courtiel's $\Delta$-activities can be generalized to polymatroids in a straightforward way, moreover, they induce a Crapo decomposition also in this case.

We show that, however, embedding activities of hypergraphs are not all $\Delta$-activities, hence Courtiel's conjecture does not hold for polymatroids, i.e. there exist activities for polymatroids that induce a Crapo-decomposition, but that are not $\Delta$-activities.

Also, we can slightly modify a hypergraph example to obtain a counterexample for Courtiel's original conjecture, that is, we give a graph and  activities that cannot be realized as $\Delta$-activities, but they yield a Crapo-decomposition.

\subsection{The definition of $\Delta$-activities}

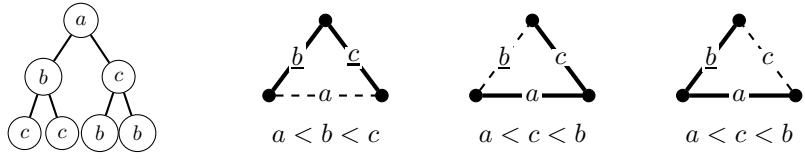
\begin{figure}
\begin{tikzpicture}[scale=.25]		
	\tikzstyle{o}=[circle,draw,scale=0.8]
	\tikzstyle{c}=[circle,draw,fill,scale=0.5]
	
	\begin{scope}[shift={(-12,0)}]
		\node [o] (1) at (4, 6) {$a$};
		\node [o] (2) at (2, 3) {$b$};
		\node [o] (3) at (6, 3) {$c$};
		\node [o] (4) at (1, 0) {$c$};
		\node [o] (5) at (3, 0) {$c$};
		\node [o] (6) at (5, 0) {$b$};
		\node [o] (7) at (7, 0) {$b$};
		
		\draw [thick] (1) -- (2);
		\draw [thick] (1) -- (3);
		\draw [thick] (2) -- (4);
		\draw [thick] (2) -- (5);
		\draw [thick] (3) -- (6);
		\draw [thick] (3) -- (7);
	\end{scope}
\begin{scope}[shift={(2,0)}]
\node [c] (1) at (0, 2) {};
\node [c] (2) at (6, 2) {};
\node [c] (3) at (3, 6) {};

\draw [dashed, thick] (1) to node[fill = white,inner sep=1pt]{$a$} (2);
\draw [ultra thick] (1) to node[fill = white,inner sep=1pt]{$\underline{b}$} (3);
\draw [ultra thick] (2) to node[fill = white,inner sep=1pt]{$\underline{c}$} (3);
\node at (3, 0) {$a < b < c$};	
\end{scope}
\begin{scope}[shift={(13,0)}]
\node [c] (1) at (0, 2) {};
\node [c] (2) at (6, 2) {};
\node [c] (3) at (3, 6) {};

\draw [ultra thick] (1) to node[fill = white,inner sep=1pt]{$a$} (2);
\draw [dashed, thick] (1) to node[fill = white,inner sep=1pt]{$\underline{b}$} (3);
\draw [ultra thick] (2) to node[fill = white,inner sep=1pt]{$c$} (3);
\node at (3, 0) {$a < c < b$};	
\end{scope}
\begin{scope}[shift={(24,0)}]
\node [c] (1) at (0, 2) {};
\node [c] (2) at (6, 2) {};
\node [c] (3) at (3, 6) {};

\draw [ultra thick] (1) to node[fill = white,inner sep=1pt]{$a$} (2);
\draw [ultra thick] (1) to node[fill = white,inner sep=1pt]{$\underline{b}$} (3);
\draw [dashed, thick] (2) to node[fill = white,inner sep=1pt]{$c$} (3);
\node at (3, 0) {$a < c < b$};
\end{scope}
\end{tikzpicture}
\caption{An example for $\Delta$-activities. Underneath each spanning tree, we have its corresponding ordering. Nontrivially active elements are underlined.}\label{fig:Delta_activities}
\end{figure}

Courtiel defines $\Delta$-activities the following way.
Take a matroid $M=(E,\mathcal{B})$. Let $\Delta$ be a rooted binary decision tree whose vertices are labeled by elements of $E$, moreover, each branch from the root to a leaf has exactly $|E|$ vertices, and each $e\in E$ occurs once as a label on the branch. (See Figure \ref{fig:Delta_activities}.) 

Then, based on $\Delta$, the following activity notion can be defined \cite{Courtiel}. Let $B\in\mathcal{B}$ be a basis. We associate an ordering $\leq_B$ to $B$ the following way. Take the element labelling the root node. If that element is in $B$, then take the left child of the root, if it is not in $B$, then take the right child of the root, and continue similarly.  The ordering $\leq_B$ is defined as the order in which we see the labels. As on each branch, each vertex appears once as a label, this is a well-defined ordering of $E$.
An element $e$ is then defined \emph{internally active} (\emph{externally active}) in $B$, if there exist no $f >_B e$ such that $B-e+f$ ($B+e-f$) is a basis. See Figure \ref{fig:Delta_activities} for an example.

\begin{remark}
	Note the important difference that while for embedding activities, one defined an element to be active if if could not be switched to a \emph{smaller} element, here, we define an edge active if it cannot be switched to a \emph{larger} element. For activities with respect to a fixed edge ordering, such a difference does not matter as one can simply reverse the ordering. However, this cannot be done for embedding or $\Delta$-activities, hence this difference has relevance here.
\end{remark}

We call an element $e$ nontrivially internally $\Delta$-active for the basis $B$ if $e\in B$ and $e$ is internally $\Delta$-active for $B$. Similarly, we call an element $e$ nontrivially externally $\Delta$-active for $B$ if $e\notin B$ and $e$ is externally $\Delta$-active for $B$. We call an element nontrivially active for $B$ if it is either nontrivially internally active for $B$ or nontrivially externally active for $B$.

We note that Courtiel actually defines $\Delta$-activities to be what we called nontrivial $\Delta$-activities above. However, as each basis has the same number of elements, considering activities instead of nontrivial activities amounts to a trivial shift in the statistics. We choose to use the above definition for activities, because for polymatroids, there is no well-defined sense of ``being in the basis'', hence the ``usual'' definition would be hard to generalize.

$\Delta$-activity can be naturally generalized for polymatroids the following way. Let $P$ be a polymatroid on ground set $E$. Let us now suppose that $P\subset \mathbb{Z}_{\geq 0}^E$, and for $e\in E$, let $r(e)=\max\{b(e): b\in P\}$.
Let $\Delta$ be a rooted tree such that each branch from the root to a leaf has exactly $|E|$ vertices and each $e\in E$ occurs exactly once as a label. Furthermore, if a vertex is labeled by $e$, then it has $r(e)+1$ children.

We associate an ordering $<_b$ to a basis $b\in P$ the following way. Let $x\in E$ be the element labeling the root of $\Delta$. If $b(x)=i$, choose the $(i+1)^{th}$ child of the root, and continue with the element labeling that vertex. (Note that $b(x)\leq r(x)$, hence this makes sense.) $<_b$ is defined as the ordering of $E$ in which we see the elements on the branch corresponding to $b$. An element $e\in E$ is internally active (resp. externally active) in $b$ if there exist no $f >_b e$ such that $b-\mathbf{1}_e+\mathbf{1}_f$ (resp. $b+\mathbf{1}_e-\mathbf{1}_f$) is a basis.

\begin{prop}
	For a polymatroid $P$ and any choice of decision tree $\Delta$,  
	\begin{enumerate}
		\item the Crapo intervals $\{C_{\Delta}(b) \mid b\in P\}$ partition $\mathbb{Z}^E$,	
		
		\item for each $c\in\mathbb{Z}^E$, if $c\in C_{\Delta}(b)$, then $d_1(P,c)=d_1(b,c)$.
	\end{enumerate}
\end{prop}
\begin{proof}
	Let us denote the Crapo interval associated to a basis $b$ using $\Delta$-activities by $C_\Delta(b)$.
	
	First we need to show that if $b_1$ and $b_2$ are distinct bases, then $C_\Delta(b_1)\cap C_\Delta(b_2)=\emptyset$. This can be proved similarly as in \cite[Theorem 9.3]{BKP}. In particular,
	take the branch of $\Delta$ corresponding to $b_1$ and the branch corresponding to $b_2$. Suppose that the last common vertex of the two brances has label $e$. Then $b_1(e)\neq b_2(e)$. By symmetry, we may assume $b_1(e)> b_2(e)$. Then, by the definition of a polymatroid, there exist $f\in E$ such that $b_1-\mathbf{1}_{e}+\mathbf{1}_{f}$ and $b_2+\mathbf{1}_{e}-\mathbf{1}_{f}$ are both bases of $P$. As the branches of $b_1$ and $b_2$ agree until the vertex labeled by $e$, the elements preceeding $e$ in $<_{b_1}$ and in $<_{b_2}$ agree, and for these elements $g$, we have $b_1(g)=b_2(g)$. Hence $f >_{b_1} e$ and $f >_{b_2} e$, thus, $e$ is internally $\Delta$-passive in $b_1$ and externally $\Delta$-passive in $b_2$. Thus, for $x\in C_\Delta(b_1)$ and $y\in C_\Delta(b_2)$ we have a $x(e) \geq b_1(e) > b_2(e) \geq y(e)$, which means $x\neq y$.
	
	Next, we show that for each $c\in \mathbb{Z}^n$, there exist a basis $b$ such that $c\in C_\Delta(b)$.
	
	For this part, we will copy the proof of Theorem \ref{thm:each_point_in_some_interval}, but this time it will be much easier.
	Let us introduce $B_c$, where $B_c=\{b \in P, d_1(P,c)=d_1(b,c)\}$. 
	Notice that in the proof of Proposition \ref{prop:external_internal_partition_for_c} we only used that hypertrees are the bases of a polymatroid. Hence for an arbitrary polymatroid, the groundset $E$ can be partitioned into elements that are external for $c$ (meaning that $b(e)\leq c(e)$ for each basis $b\in B_c$), elements that are internal for $c$ (meaning that $b(e) \geq c(e)$ for each basis $b\in B_c$), and elements where $b(e)=c(e)$ for each $b\in B_c$.
	
	We define an ordering of the bases in $B_c$ the following way.
	Let $b, b'\in B_c$, and suppose that the branches of $\Delta$ corresponding to them part at a vertex labeled $e$. Then we have $b(e)\neq b'(e)$. We define $b' <_c b$ if
	\begin{enumerate}
		\item $e$ is external for $c$ and $b'(e) > b(e)$,
		\item $e$ is internal for $c$ and $b'(e) < b(e)$.
	\end{enumerate}
    Note that $<_c$ is a complete ordering on $B_c$.
    
    Take an arbitrary basis $b\in B_c$. We show that if $c\notin C_\Delta(b)$, then we can find some $b'\in B_c$ with $b' <_c b$. As $B_c$ is finite, this means that we eventually find a basis $b''\in B_c$ with $c\in C_\Delta(b'')$.
    
    Suppose that $c\notin C_\Delta(b)$. This means that either there are some element(s) $e$ such that $e$ is externally $\Delta$-passive in $b$ but $c(e)>b(e)$, or some element(s) $e$ such that $e$ is internally $\Delta$-passive in $b$ but $c(e)<b(e)$ (or both).
    Take the first such bad element $e$ according to $<_b$. 
    
    Suppose that $e$ is externally $\Delta$-passive in $b$ but $c(e)>b(e)$. Then by the definition of exernal passivity, there exist some $e'>_b e$ such that $b'=b+\mathbf{1}_e-\mathbf{1}_{e'}$ is also a basis. Since $b'(e)> b(e)$ and $c(e)>b(e)$, we have $d_1(b',c)\leq d_1(b,c)$, hence $b'\in B_c$. As the branch of $b$ contains the label $e$ before label $e'$, the branches of $b$ and $b'$ part at $e$. We conclude that $b' <_c b$.
    
    The case when the first bad element $e$ is internal is completely analogous.
\end{proof}

Courtiel shows that for graphs, embedding activities fit into the framework of $\Delta$-activities. However, for hypergraphs we show an example that embedding activities are not all $\Delta$-activities.

We will use the following observation, whose matroid version was noted by Courtiel \cite{Courtiel}. For a polymatroid, we call an element $e$ nontrivially internally active for a basis $b\in P$, if there exist a basis $b'\in P$ with $b'(e) < b(e)$ and $e$ is internally active for $b$. We call $e$ nontrivially externally active for $b$, if there exist a basis $b'\in P$ with $b'(e)>b(e)$ 
and $e$ is externally active for $b$. 

\begin{prop}\label{prop:Delta-necessary_condition}
	For any polymatroid $\Delta$-activity, there exist an element of the ground set that is neither nontrivially internally $\Delta$-active, nor nontrivially externally $\Delta$-active for any basis.
\end{prop}
\begin{proof}
The element $e$ in the root of the decision tree will be the first element for each ordering $<_B$. It follows directly from the definition of a polymatroid that this element is never nontrivially internally active, nor nontrivially externally active for any basis.
\end{proof}

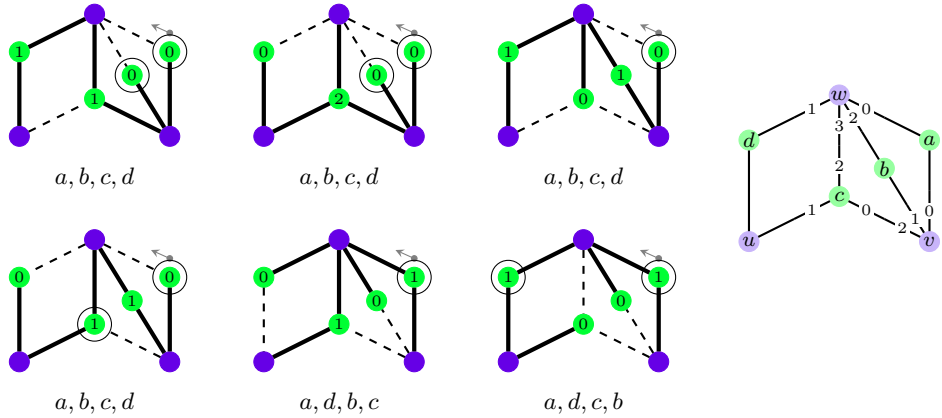
\begin{figure}[t]
	\begin{tikzpicture}[scale=.25]		
		\tikzstyle{o}=[circle,draw,scale=0.8]
		
		\begin{scope}[shift={(-13,0)}]
			\node [o,fill=e,e] (a) at (22, 6.5) {};
			\node [o,fill=e,e] (c) at (18, 4) {};
			\node [o,fill=e,e] (d) at (14, 6.5) {};
			\node [o,fill=e,e] (b) at (20,5.25) {};
			\node [o,fill=v,v] (u) at (14, 2) {};
			\node [o,fill=v,v] (v) at (22, 2) {};
			\node [o,fill=v,v] (w) at (18, 8.5) {};
			
			\draw [ultra thick] (w) -- (d);
			\draw [ultra thick] (d) -- (u);
			\draw [dashed, thick] (w) -- (b);
			\draw [ultra thick] (w) -- (c);
			\draw [dashed, thick] (u) -- (c);
			\draw [dashed, thick] (w) -- (a);
			\draw [ultra thick] (a) -- (v);
			\draw [ultra thick] (v) -- (c);
			\draw [ultra thick] (v) -- (b);
			
			\draw [fill,color=gray] (22,7.5) circle [radius=.15];
			\draw [->,>=stealth,color=gray] (22,7.5) -- (21,7.95);
			
			\node at (22, 6.5) {{\tiny{$0$}}};
			\node at (18, 4) {{\tiny{$1$}}};
			\node at (20,5.25) {{\tiny{$0$}}};
			\node at (14, 6.5) {{\tiny{$1$}}};
			\draw (22, 6.5) circle [radius=0.9];
	        \draw (20, 5.25) circle [radius=0.9];
	        		
			\node at (18, -0.2) {{\small{$a,b,c,d$}}};
		\end{scope}
		
		\begin{scope}[shift={(0,0)}]
			\node [o,fill=e,e] (a) at (22, 6.5) {};
			\node [o,fill=e,e] (c) at (18, 4) {};
			\node [o,fill=e,e] (d) at (14, 6.5) {};
			\node [o,fill=e,e] (b) at (20,5.25) {};
			\node [o,fill=v,v] (u) at (14, 2) {};
			\node [o,fill=v,v] (v) at (22, 2) {};
			\node [o,fill=v,v] (w) at (18, 8.5) {};
			
			\draw [dashed, thick] (w) -- (d);
			\draw [ultra thick] (d) -- (u);
			\draw [dashed, thick] (w) -- (b);
			\draw [ultra thick] (w) -- (c);
			\draw [ultra thick] (u) -- (c);
			\draw [dashed, thick] (w) -- (a);
			\draw [ultra thick] (a) -- (v);
			\draw [ultra thick] (v) -- (c);
			\draw [ultra thick] (v) -- (b);
			
			\draw [fill,color=gray] (22,7.5) circle [radius=.15];
			\draw [->,>=stealth,color=gray] (22,7.5) -- (21,7.95);
			
			\node at (22, 6.5) {{\tiny{$0$}}};
			\node at (18, 4) {{\tiny{$2$}}};
			\node at (20,5.25) {{\tiny{$0$}}};
			\node at (14, 6.5) {{\tiny{$0$}}};
			\draw (22, 6.5) circle [radius=0.9];
			\draw (20, 5.25) circle [radius=0.9];
			
			\node at (18, -0.2) {{\small{$a,b,c,d$}}};
		\end{scope}
		
		\begin{scope}[shift={(13, 0)}]
			\node [o,fill=e,e] (a) at (22, 6.5) {};
			\node [o,fill=e,e] (c) at (18, 4) {};
			\node [o,fill=e,e] (d) at (14, 6.5) {};
			\node [o,fill=e,e] (b) at (20,5.25) {};
			\node [o,fill=v,v] (u) at (14, 2) {};
			\node [o,fill=v,v] (v) at (22, 2) {};
			\node [o,fill=v,v] (w) at (18, 8.5) {};
			
		    \draw [ultra thick] (w) -- (d);
		    \draw [ultra thick] (d) -- (u);
		    \draw [ultra thick] (w) -- (b);
		    \draw [ultra thick] (w) -- (c);
		    \draw [dashed, thick] (u) -- (c);
		    \draw [dashed, thick] (w) -- (a);
			\draw [ultra thick] (a) -- (v);
			\draw [dashed, thick] (v) -- (c);
			\draw [ultra thick] (v) -- (b);
			
			\draw [fill,color=gray] (22,7.5) circle [radius=.15];
			\draw [->,>=stealth,color=gray] (22,7.5) -- (21,7.95);
			
			\node at (22, 6.5) {{\tiny{$0$}}};
			\node at (18, 4) {{\tiny{$0$}}};
			\node at (20,5.25) {{\tiny{$1$}}};
			\node at (14, 6.5) {{\tiny{$1$}}};
			\draw (22, 6.5) circle [radius=0.9];
			
			\node at (18, -0.2) {{\small{$a,b,c,d$}}};
		\end{scope}
		
		\begin{scope}[shift={(-13, -12)}]
			\node [o,fill=e,e] (a) at (22, 6.5) {};
			\node [o,fill=e,e] (c) at (18, 4) {};
			\node [o,fill=e,e] (d) at (14, 6.5) {};
			\node [o,fill=e,e] (b) at (20,5.25) {};
			\node [o,fill=v,v] (u) at (14, 2) {};
			\node [o,fill=v,v] (v) at (22, 2) {};
			\node [o,fill=v,v] (w) at (18, 8.5) {};
			
			\draw [dashed, thick] (w) -- (d);
			\draw [ultra thick] (d) -- (u);
			\draw [ultra thick] (w) -- (b);
			\draw [ultra thick] (w) -- (c);
			\draw [ultra thick] (u) -- (c);
			\draw [dashed, thick] (w) -- (a);
			\draw [ultra thick] (a) -- (v);
			\draw [dashed, thick] (v) -- (c);
			\draw [ultra thick] (v) -- (b);
			
			\draw [fill,color=gray] (22,7.5) circle [radius=.15];
			\draw [->,>=stealth,color=gray] (22,7.5) -- (21,7.95);
			
			\node at (22, 6.5) {{\tiny{$0$}}};
			\node at (18, 4) {{\tiny{$1$}}};
			\node at (20,5.25) {{\tiny{$1$}}};
			\node at (14, 6.5) {{\tiny{$0$}}};
			\draw (22, 6.5) circle [radius=0.9];
			\draw (18, 4) circle [radius=0.9];
			
			\node at (18, -0.2) {{\small{$a,b,c,d$}}};
		\end{scope}
		
		\begin{scope}[shift={(0, -12)}]
			\node [o,fill=e,e] (a) at (22, 6.5) {};
			\node [o,fill=e,e] (c) at (18, 4) {};
			\node [o,fill=e,e] (d) at (14, 6.5) {};
			\node [o,fill=e,e] (b) at (20,5.25) {};
			\node [o,fill=v,v] (u) at (14, 2) {};
			\node [o,fill=v,v] (v) at (22, 2) {};
			\node [o,fill=v,v] (w) at (18, 8.5) {};
			
			\draw [ultra thick] (w) -- (d);
			\draw [dashed, thick] (d) -- (u);
			\draw [ultra thick] (w) -- (b);
			\draw [ultra thick] (w) -- (c);
			\draw [ultra thick] (u) -- (c);
			\draw [ultra thick] (w) -- (a);
			\draw [ultra thick] (a) -- (v);
			\draw [dashed, thick] (v) -- (c);
			\draw [dashed, thick] (v) -- (b);
			
			\draw [fill,color=gray] (22,7.5) circle [radius=.15];
			\draw [->,>=stealth,color=gray] (22,7.5) -- (21,7.95);
			
			\node at (22, 6.5) {{\tiny{$1$}}};
			\node at (18, 4) {{\tiny{$1$}}};
			\node at (20,5.25) {{\tiny{$0$}}};
			\node at (14, 6.5) {{\tiny{$0$}}};
			\draw (22, 6.5) circle [radius=0.9];
			
			\node at (18, -0.2) {{\small{$a,d,b,c$}}};
		\end{scope}
		
		\begin{scope}[shift={(13, -12)}]
			\node [o,fill=e,e] (a) at (22, 6.5) {};
			\node [o,fill=e,e] (c) at (18, 4) {};
			\node [o,fill=e,e] (d) at (14, 6.5) {};
			\node [o,fill=e,e] (b) at (20,5.25) {};
			\node [o,fill=v,v] (u) at (14, 2) {};
			\node [o,fill=v,v] (v) at (22, 2) {};
			\node [o,fill=v,v] (w) at (18, 8.5) {};
			
			\draw [ultra thick] (w) -- (d);
			\draw [ultra thick] (d) -- (u);
			\draw [ultra thick] (w) -- (b);
			\draw [dashed, thick] (w) -- (c);
			\draw [ultra thick] (u) -- (c);
			\draw [ultra thick] (w) -- (a);
			\draw [ultra thick] (a) -- (v);
			\draw [dashed, thick] (v) -- (c);
			\draw [dashed, thick] (v) -- (b);
			
			\draw [fill,color=gray] (22,7.5) circle [radius=.15];
			\draw [->,>=stealth,color=gray] (22,7.5) -- (21,7.95);
			
			\node at (22, 6.5) {{\tiny{$1$}}};
			\node at (18, 4) {{\tiny{$0$}}};
			\node at (20,5.25) {{\tiny{$0$}}};
			\node at (14, 6.5) {{\tiny{$1$}}};
			\draw (22, 6.5) circle [radius=0.9];
			\draw (14, 6.5) circle [radius=0.9];
			
			\node at (18, -0.2) {{\small{$a,d,c,b$}}};
		\end{scope}
		\begin{scope}[shift={(23, -6)}, scale=1.2]
			\node [o,fill=le,le] (a) at (22, 6.5) {};
			\node [o,fill=le,le] (c) at (18, 4) {};
			\node [o,fill=le,le] (d) at (14, 6.5) {};
			\node [o,fill=le,le] (b) at (20,5.25) {};
			\node [o,fill=lv,lv] (u) at (14, 2) {};
			\node [o,fill=lv,lv] (v) at (22, 2) {};
			\node [o,fill=lv,lv] (w) at (18, 8.5) {};
			
			\draw [thick] (w) to node[fill = white,inner sep=1pt]{\tiny $1$} (16,7.5);
			\draw [thick] (16, 7.5) -- (d);
			\draw [thick] (d) -- (u);
			\draw [thick] (w) to node[fill = white,inner sep=1pt]{\tiny $2$} (19, 6.875);
			\draw [thick] (19, 6.875) -- (b);
			\draw [thick] (w) to node[fill = white,inner sep=1pt]{\tiny $3$} (18, 6.25); 
			\draw [thick] (u) -- (16,3);
			\draw [thick] (16,3) to node[fill = white,inner sep=1pt]{\tiny $1$} (c); 
			\draw [thick] (18, 6.25) to node[fill = white,inner sep=1pt]{\tiny $2$} (c);
			\draw [thick] (w) to node[fill = white,inner sep=1pt]{\tiny $0$} (20, 7.5);
			\draw [thick] (20, 7.5) -- (a);
			\draw [thick] (a) -- (22,4.25);
			\draw [thick] (22, 4.25) to node[fill = white,inner sep=1pt]{\tiny $0$} (v);
			\draw [thick] (v) to node[fill = white,inner sep=1pt]{\tiny $2$} (20,3);
			\draw [thick] (20,3) to node[fill = white,inner sep=1pt]{\tiny $0$} (c);
			\draw [thick] (b) -- (21,3.625);
			\draw [thick] (21, 3.625) to node[fill = white,inner sep=1pt]{\tiny $1$} (v);
			
			\node at (22, 6.5) {{\footnotesize{$a$}}};
			\node at (18, 4) {{\footnotesize{$c$}}};
			\node at (20,5.25) {{\footnotesize{$b$}}};
			\node at (14, 6.5) {{\footnotesize{$d$}}};
			
			\node at (14, 2) {{\footnotesize{$u$}}};
			\node at (22, 2) {{\footnotesize{$v$}}};
			\node at (18, 8.5) {{\footnotesize{$w$}}};
		\end{scope}
	\end{tikzpicture}
	\caption{A counterexample for the polymatroid version of Courtiel's conjecture. See the proof of Proposition \ref{prop:hypergraph_Courtiel_counterex} for the notations and more explanation.}\label{fig:Delta_counterex}
\end{figure}

\begin{prop}\label{prop:hypergraph_Courtiel_counterex}
	Embedding activities of hypergraphs are not all $\Delta$-activities. Hence there exist 
	activities for polymatroids yielding a Crapo decomposition that are not $\Delta$-activities. 
\end{prop}
\begin{proof}
	Take the hypergraph represented by the bipartite graph on Figure \ref{fig:Delta_counterex}, with 3 vertices (violet nodes) and 4 hyperedges (emerald nodes). The labelings of the nodes and the ribbon structure are depicted on the rightmost panel. (The cyclic ordering of the edges incident to a node is indicated by the numbers written on the edges near the node. For clarity of the picture, the (trivial) information is omitted for degree 2 nodes.)
	Let the basis be $(a,aw)$. 
	
	The 6 left panels of Figure \ref{fig:Delta_counterex} show the 6 hypertrees of the hypergraph, together with their representing Jaeger trees, and the corresponding orderings $\leq_h$ on the hyperedges (that is, emerald nodes). For each hypertree, nontrivially (internally or externally) embedding active emerald nodes are circled.

    By Proposition \ref{prop:Delta-necessary_condition}, for a $\Delta$-activity, there must exist a hyperedge that is neither nontrivially internally active, nor nontrivially externally active for any basis.
	
	However, we see on Figure \ref{fig:Delta_counterex} that for the constructed example, each emerald node $a,b,c$ and $d$ is nontrivially embedding-active for some hypertree, hence we conclude that no decision tree can give us this set of active hyperedges.
\end{proof}

We can use the example of Figure \ref{fig:Delta_counterex} to also construct a counterexample for Courtiel's conjecture on graphs.

\begin{figure}[t]
	\begin{tikzpicture}[scale=.2]		
		\tikzstyle{o}=[circle,draw,scale=0.7]
		
		\begin{scope}[shift={(-24,0)}]
			\node [o,fill] (u) at (0, 2) {};
			\node [o,fill] (v) at (8, 2) {};
			\node [o,fill] (w) at (4, 8.5) {};
			
			\draw [dashed, thick, bend right = 20] (v)  to node[fill = white,inner sep=1pt]{\underline{$a$}} (w);
			\draw [dashed, thick, bend left = 20] (v) to node[fill = white,inner sep=1pt]{\underline{$b$}} (w);
			\draw [ultra thick] (u) to node[fill = white,inner sep=1pt]{$c$} (v);
			\draw [ultra thick] (u) to node[fill = white,inner sep=1pt]{$d$} (w);
						
			\node at (4, -0.2) {{\small{$a<b<c<d$}}};
		\end{scope}
		
		\begin{scope}[shift={(-12,0)}]
			\node [o,fill] (u) at (0, 2) {};
			\node [o,fill] (v) at (8, 2) {};
			\node [o,fill] (w) at (4, 8.5) {};
			
			\draw [dashed, thick, bend right = 20] (v)  to node[fill = white,inner sep=1pt]{\underline{$a$}} (w);
			\draw [ultra thick, bend left = 20] (v) to node[fill = white,inner sep=1pt]{$b$} (w);
			\draw [dashed, thick] (u) to node[fill = white,inner sep=1pt]{$c$} (v);
			\draw [ultra thick] (u) to node[fill = white,inner sep=1pt]{$d$} (w);
			
			\node at (4, -0.2) {{\small{$a<b<c<d$}}};
		\end{scope}
		
		\begin{scope}[shift={(0, 0)}]
			
			\node [o,fill] (u) at (0, 2) {};
			\node [o,fill] (v) at (8, 2) {};
			\node [o,fill] (w) at (4, 8.5) {};
			
			\draw [dashed, thick, bend right = 20] (v)  to node[fill = white,inner sep=1pt]{\underline{$a$}} (w);
			\draw [ultra thick, bend left = 20] (v) to node[fill = white,inner sep=1pt]{$b$} (w);
			\draw [ultra thick] (u) to node[fill = white,inner sep=1pt]{\underline{$c$}} (v);
			\draw [dashed, thick] (u) to node[fill = white,inner sep=1pt]{$d$} (w);
			
			\node at (4, -0.2) {{\small{$a<b<c<d$}}};
		\end{scope}
		
		\begin{scope}[shift={(12, 0)}]
			
			\node [o,fill] (u) at (0, 2) {};
			\node [o,fill] (v) at (8, 2) {};
			\node [o,fill] (w) at (4, 8.5) {};
			
			\draw [ultra thick, bend right = 20] (v)  to node[fill = white,inner sep=1pt]{\underline{$a$}} (w);
			\draw [dashed, thick, bend left = 20] (v) to node[fill = white,inner sep=1pt]{$b$} (w);
			\draw [ultra thick] (u) to node[fill = white,inner sep=1pt]{$c$} (v);
			\draw [dashed, thick] (u) to node[fill = white,inner sep=1pt]{$d$} (w);
			
			\node at (4, -0.2) {{\small{$a<d<b<c$}}};
		\end{scope}
		
		\begin{scope}[shift={(24, 0)}]
			\node [o,fill] (u) at (0, 2) {};
			\node [o,fill] (v) at (8, 2) {};
			\node [o,fill] (w) at (4, 8.5) {};
			
			\draw [ultra thick, bend right = 20] (v)  to node[fill = white,inner sep=1pt]{\underline{$a$}} (w);
			\draw [dashed, thick, bend left = 20] (v) to node[fill = white,inner sep=1pt]{$b$} (w);
			\draw [dashed, thick] (u) to node[fill = white,inner sep=1pt]{$c$} (v);
			\draw [ultra thick] (u) to node[fill = white,inner sep=1pt]{\underline{$d$}} (w);
			
			\node at (4, -0.2) {{\small{$a<d<c<b$}}};
		\end{scope}
	\end{tikzpicture}
	\caption{A counterexample for Courtiel's conjecture. Activities are defined with respect to each tree's associated ordering. The nontrivially active edges are underlined.
		See the proof of Proposition \ref{prop:Courtiel_counterex_graph} for more explanation.}\label{fig:graph_Delta_counterex}
\end{figure}
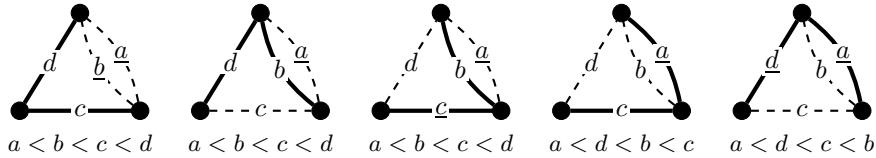

\begin{prop}\label{prop:Courtiel_counterex_graph}
	There exist a graph with a notion of activity that induces a Crapo-decomposition, but that is not a $\Delta$-activity.
\end{prop}
\begin{proof}
	Take only those hypertrees of Figure \ref{fig:Delta_counterex} that are $0/1$-valued. These hypertrees are known to form a matroid in general (the hypergraphic matroid), but in this special case, they even form a graphic matroid. Indeed, take the graph on Figure \ref{fig:graph_Delta_counterex}. Its five spanning trees are $cd, bd, bc, ac$ and $ad$, corresponding to the five $0/1$-valued hypertrees of Figure \ref{fig:Delta_counterex}. Let us assign to these spanning trees the orders seen on Figure \ref{fig:Delta_counterex}. That is, the orders corresponding to $cd, bd$ and $bc$ are $a<b<c<d$, the order correspondng to $ac$ is $a<d<b<c$ and the order corresponding to $ad$ is $a<d<c<b$. Let's define activities in each tree with respect to the tree's associated ordering, using the minimum rule (that is, an element is internally active if its value in the basis cannot be decreased, and increased on a smaller element instead). 
	On Figure \ref{fig:graph_Delta_counterex}, we underlined the nontrivially (internally or externally) active edges for each tree.
	
	It is easy to check by hand that these activities yield a Crapo-decomposition. For the nontrivial part, see the following table, checking that each vector of $\{0,1\}^E$ is in the Crapo-interval of exactly one spanning tree. The vectors corresponding to spanning trees are bold, and the nontrivially (internally or externally) active coordinates of the bases are underlined. After each vector, you can see which tree's interval contains it.
	
	It is also easy to check that for each vector, the (characteristic vector of the) spanning tree containing it in its Crapo-interval has minimal $d_1$ distance to the vector among the spanning trees.
	\medskip
	
	\begin{tabular}{cllcllcllcl}
		0000 & {\color{blue} ad \quad} & & 0100 & {\color{dy} bc \quad} & & 1000 & {\color{blue} ad \quad} & & 1100 & {\color{dy} bc} \\
		0001 & {\color{blue} ad} & & \textbf{\underline{0}101} & bd & & \textbf{\underline{1}00\underline{1}} & {\color{blue} ad} & & 1101 & bd \\
		0010 & {\color{red}ac} & & \textbf{\underline{0}1\underline{1}0} & {\color{dy} bc} & & \textbf{\underline{1}010} & {\color{red}ac} & & 1110 & {\color{dy} bc} \\
		\textbf{\underline{0}\underline{0}11} & {\color{de} cd} & & 0111 & {\color{de} cd} & & 1011 & {\color{de} cd} & & 1111 & {\color{de} cd}
	\end{tabular}

\noindent
We can see that it is still true that each edge is nontrivially active for some spanning tree, hence by Proposition \ref{prop:Delta-necessary_condition}, this activity notion cannot be obtained as a $\Delta$-activity.
\end{proof}

\begin{remark}
	One could construct a notion of activities for each hypergraphic polymatroid mimicking the proof of Propositon \ref{prop:Courtiel_counterex_graph}. (That is, for each $0/1$-valued hypertree $h$, defining activities for $h$ using the ordering $<_h$.) These will not always induce a Crapo-decomposition.
\end{remark}

\bibliographystyle{plain}
\bibliography{Bernardi}

\begin{thebibliography}{10}

\bibitem{Bernardi_first}
Olivier Bernardi.
\newblock A characterization of the {T}utte polynomial via combinatorial
  embeddings.
\newblock {\em Ann. Comb.}, 12(2):139--153, 2008.

\bibitem{BKP}
Olivier Bernardi, Tam\'{a}s K\'{a}lm\'{a}n, and Alexander Postnikov.
\newblock Universal {T}utte polynomial.
\newblock {\em Adv. Math.}, 402:Paper No. 108355, 74, 2022.

\bibitem{Cameron-Fink}
A.~Cameron and A.~Fink.
\newblock A lattice point counting generalisation of the tutte polynomial.
\newblock https://arxiv.org/pdf/1604.00962.pdf, 2016.

\bibitem{Courtiel}
J.~Courtiel.
\newblock A general notion of activity for the tutte polynomial.
\newblock arXiv:1412.2081, 2014.

\bibitem{Edmonds_polymatroid}
Jack Edmonds.
\newblock Submodular functions, matroids, and certain polyhedra [mr0270945].
\newblock In {\em Combinatorial optimization---{E}ureka, you shrink!}, volume
  2570 of {\em Lecture Notes in Comput. Sci.}, pages 11--26. Springer, Berlin,
  2003.

\bibitem{hiperTutte}
Tam\'{a}s K\'{a}lm\'{a}n.
\newblock A version of {T}utte's polynomial for hypergraphs.
\newblock {\em Adv. Math.}, 244:823--873, 2013.

\bibitem{hyperBernardi}
Tam\'{a}s K\'{a}lm\'{a}n and Lilla T\'othm\'er\'esz.
\newblock Hypergraph polynomials and the {B}ernardi process.
\newblock {\em Algebraic Combinatorics}, 3(5):1099--1139, 2020.

\bibitem{semibalanced}
Tam\'{a}s K\'{a}lm\'{a}n and Lilla T\'{o}thm\'{e}r\'{e}sz.
\newblock Root polytopes and {J}aeger-type dissections for directed graphs.
\newblock {\em Mathematika}, 68(4):1176--1220, 2022.

\bibitem{Schrijver_CombOpt}
Alexander Schrijver.
\newblock {\em Combinatorial optimization. {P}olyhedra and efficiency. {V}ol.
  {B}}, volume~24 of {\em Algorithms and Combinatorics}.
\newblock Springer-Verlag, Berlin, 2003.
\newblock Matroids, trees, stable sets, Chapters 39--69.

\bibitem{Tutte}
W.~T. Tutte.
\newblock A contribution to the theory of chromatic polynomials.
\newblock {\em Canad. J. Math.}, 6:80--91, 1954.

\end{thebibliography}

\end{document}